\documentclass[12pt]{article}
\usepackage{amssymb,graphicx,subfigure}
\usepackage{amsmath,amsfonts,enumerate,amsthm}
\usepackage{subfigure}
\usepackage[letterpaper,margin=1in]{geometry}

\def\eps{\varepsilon}

\begin{document}

\title{Computing Slow Manifolds of Saddle Type}
\author{John Guckenheimer
\thanks{Mathematics Department, Cornell University, Ithaca, NY 14853}
\and
Christian Kuehn
\thanks{Center for Applied Mathematics, Cornell University, Ithaca, NY 14853}}

\maketitle

\begin{abstract}
Slow manifolds are important geometric structures in the state spaces of 
dynamical systems with multiple time scales. This paper introduces an 
algorithm for computing trajectories on slow manifolds that are normally
hyperbolic with both stable and unstable fast manifolds. We present two 
examples of bifurcation problems where these manifolds play a key role 
and a third example in which saddle-type slow manifolds are part of a
traveling wave profile of a partial differential equation. Initial value 
solvers are incapable of computing trajectories on saddle-type slow manifolds, 
so the slow manifold of saddle type (SMST) algorithm presented here 
is formulated as a boundary value method. 
We take an empirical approach here to assessing the accuracy and effectiveness
of the algorithm.
\end{abstract}

\section{Introduction}
\emph{Slow-fast} vector fields have the form
\begin{equation}
\begin{split}
\eps \dot{x} & =  f(x,y,\eps)  \\
\dot{y} & =  g(x,y,\eps)
\end{split}
\label{sfs}
\end{equation}
with $x \in R^m$ the fast variable, $y \in R^n$ the slow
variable and $\eps$ a small parameter that represents the ratio of
time scales. The pair $(x,y)$ will be denoted by $z$ and the vector
field will be written $\dot{z} = F(z)$.
Simulation of these systems with explicit numerical 
integration algorithms is limited to time steps that are $O(\eps)$
due to numerical instabilities. However, invariant {\em slow
manifolds} on which the motion of the system has 
speed that is $O(1)$ are a common feature of slow-fast systems. 
Indeed, trajectories often spend most of their time following 
stable slow manifolds. Implicit ``stiff'' integration 
methods~\cite{HW} compute trajectories along the stable slow manifolds,
taking time steps that are $O(1)$ while avoiding the  numerical 
instabilities of explicit methods. However, no initial value solver
will compute forward trajectories that evolve on unstable slow
manifolds because the geometric instability of these trajectories is such
that nearby initial conditions diverge from one another at exponential
rates commensurate with the fast time scale. Even an
exact initial value solver in the presence of round-off errors 
of magnitude $\delta$ will amplify this round-off error to unit 
magnitude in a time that is $O(-\eps \log (\delta))$. Trajectories on slow
manifolds that expand in all normal directions can be computed by 
reversing time, but different strategies are needed to compute 
trajectories that lie on slow manifolds of saddle type. This paper 
presents the first algorithms that directly compute accurate trajectories 
of slow manifolds of saddle type.\footnote{Continuation methods such as
AUTO that follow one parameter families of periodic or homoclinic
orbits have been used to compute trajectories with saddle canards.
Their robustness leaves much to be desired as is discussed below.}

The existence of {\em normally hyperbolic} slow manifolds is established by {\em Fenichel theory} \cite{Fenichel,Jones}. The singular limit $\eps = 0$
of system~\eqref{sfs} is a differential algebraic equation with 
trajectories confined to the {\em critical manifold} $S = S_0$ defined by $f=0$.
At points of $S$ where $D_xf$ is a regular $m \times m$ matrix, the
implicit function theorem implies that $S$ is locally the graph
of a function $x=h(y)$. This equation yields the vector field $\dot{y}
=  g(h(y),y,0)$ for the {\em slow flow} on $S$. The geometry is more
complicated at {\em fold points} of $S$ where $D_xf$ is singular. It
is often possible to extend the slow flow to the fold points 
after a rescaling of the vector field~\cite{Montreal02}. Where all eigenvalues of
$D_xf$ have nonzero real parts, Fenichel proved the existence of 
invariant slow manifolds $S_\eps$ for $\eps > 0$ small. These 
{\em normally hyperbolic} slow manifolds are within an $O(\eps)$
distance from the critical manifold $S_0$ and the flow on $S_\eps$
converges to the slow flow on $S_0$ as $\eps \to 0$. Fenichel theory
is usually developed in the context of {\em overflowing} slow
manifolds with boundaries which trajectories may leave through the
boundaries. In this setting, slow manifolds are not unique, but the
distance between a pair of slow manifolds is ``exponentially small'',
i.e. of order $O(\exp(-c/\eps))$ for a suitable positive $c$,
independent of $\eps$~\cite{Jones}.

\section{The SMST Algorithm}

This section describes a collocation method for computing slow 
manifolds of saddle type in slow-fast systems that we call the 
SMST algorithm. The numerical analysis
employed in the algorithm is straightforward; the subtlety of the problem 
appears in the formulation of discrete systems of equations with
well-conditioned Jacobians. The crucial part of the geometry is to 
specify boundary conditions for trajectory segments on a slow 
manifold that yield well-conditioned discretizations.

\begin{figure}[!htb]
  \begin{center}
    \includegraphics[height=5in]{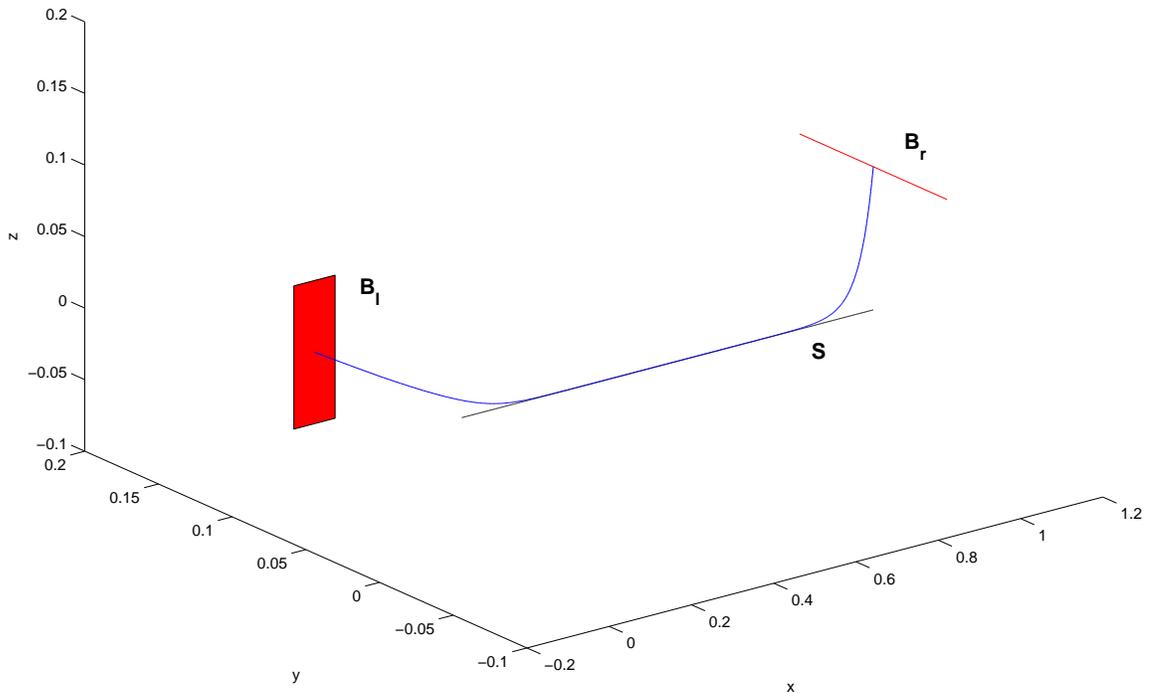}
 \end{center}
  \caption{Boundary conditions for the SMST algorithm are illustrated with 
a three dimensional example with one slow and two fast variables. The slow 
manifold of saddle type is drawn black 
and labeled $S$. A trajectory that approaches the slow manifold along a strong
stable direction and departs along a strong unstable manifold is drawn blue.
The initial point of this trajectory lies in a two dimensional manifold $B_l$ 
transverse to the stable manifold of $S$, and the final point lies in a one 
dimensional manifold $B_r$ transverse to a unstable manifold of $S$. 
}
\label{bv_su}
\end{figure}

A trajectory segment $\gamma:[a,b] \to R^{m+n}$ of system \eqref{sfs} 
is determined by its initial point $\gamma(a)$ or by another set of 
$m+n$ boundary conditions. Trajectories that follow a slow manifold 
for some distance approach the manifold initially at a fast
exponential rate and then diverge from the manifold at a fast
exponential rate. Such trajectories will be found as solutions to a
boundary value problem that imposes constraints on both $\gamma(a)$
and $\gamma(b)$.
At $\eps = 0$, there are specific arrival and
departure points. The singular limit of the trajectories we seek are
{\em candidates} $\gamma_0$ that 
consist of a fast initial segment approaching the critical manifold
$S$ along a strong stable manifold, followed by a slow segment along $S$, followed
by a fast segment that leaves $S$ along a strong unstable manifold. See Figure~\ref{bv_su}.
The initial and/or final segments may be absent.
For small $\eps > 0$, we seek $m+n$ boundary conditions that determine
a unique trajectory near the candidate. 
Initial conditions that do not lie in the strong stable manifold of a point
$p \in S$ will diverge from the slow manifold $S$ at a fast exponential
rate.  Therefore trajectories that follow the slow manifold have initial
conditions that are exponentially close to the (unknown) stable
manifold of $S$. Similarly, when trajectories depart from $S$, they 
remain exponentially close to the unstable manifold of $S$ for times
that are $O(1)$ on the fast time scale. Consequently, to have a
solvable boundary problem that is well posed, the $m+n$
boundary conditions should consist of two manifolds of dimension $k$
and $m+n-k$, the first transverse to the stable manifold of $S$ and 
the second transverse to the unstable manifold of $S$. 
Thus $u \le k\le n+u$ where $u$ is the dimension of the strong unstable 
manifolds of $S$.

The fast
segments of trajectories that precede and follow segments along the
slow manifold are readily computed with an initial value solver; the 
challenge is to locate the slow portion of the trajectory. Therefore,
the algorithm presented here takes as its input a (discretized) trajectory 
$\gamma_0:[a,b] \to C$ of the slow flow on the critical manifold
together with two submanifolds $B_l$ and $B_r$ of dimensions $k$
and $m+n-k$ that pass close to the initial and final points 
$p = \gamma_0(a)$ and $q = \gamma_0(b)$ of $\gamma_0$. The manifold
$B_l$ is assumed to be transverse to the stable manifold of $C$ and 
the manifold $B_r$ is assumed to be transverse to the unstable manifold of $C$. See Figure~\ref{bv_su}.

Denote the mesh points in the discretization of $\gamma_0$ by $a = t_0 <
t_1  < \cdots < t_N = b$. The algorithm described here is based on a
set of $(N+1)(m+n)$ equations $E$ for $\gamma(t_i)$ that yield an
approximate solution to the
boundary value problem that $\gamma$ is a trajectory of
system~\eqref{sfs} connecting $\gamma(a) \in B_l$ to
$\gamma(b) \in B_r$\footnote{If $b-a$ is allowed to vary, then the 
number of boundary conditions is increased by one.}. 
The discretized equations are based upon interpolation with cubic Hermite 
splines, though higher degree splines can be used in conjunction with
automatic differentiation methods that compute Taylor polynomials 
of the vector field at mesh points~\cite{GM}.
From points $z_j = z(t_j) \in R^{m+n}$, a $C^1$ cubic spline $\sigma$
is constructed with the $z_j$ as knot points and tangent vectors
$F(z_j)$ at these points. On the mesh interval $[t_{j-1},t_j]$,
$\sigma$ is a cubic curve whose coefficients are linear combinations
of $z_{j-1}, z_j, F(z_{j-1}), F(z_j)$ that are readily determined.
Each of the $N$ mesh intervals $[t_{j-1},t_j]$ contributes
$(m+n)$ equations to the system $E$ by requiring that 
$F(\sigma((t_{j-1}+t_j)/2) = \sigma'((t_{j-1}+t_j)/2)$. 
The values of $\sigma$ and $\sigma'$ in these equations can be expressed as
\begin{equation}
\label{spline}
\begin{split}
\sigma(\frac{t_{j-1}+t_j}{2}) & = \frac{z_{j-1}+z_j}{2} - 
\frac{(t_{j}-t_{j-1})(F(z_{j})-F(z_{j-1}))}{8} \\
\sigma'(\frac{t_{j-1}+t_j}{2}) & = \frac{3(z_{j}-z_{j-1})}{2(t_{j}-t_{j-1})}
-\frac{F(z_{j})+F(z_{j-1})}{4}
\end{split}
\end{equation}
Implicit
equations for the boundary value manifolds constitute the remaining 
$m+n$ equations in $E$. The system $E$ is solved with Newton's method
starting with the data in $\gamma_0$. The Jacobian of the system $E$
can be computed, using the derivatives of the equations~\eqref{spline} 
with respect to $z_{j-1}, z_j$.

Two types of error estimates are of interest for this algorithm. 
On each mesh interval, there is a local error estimate for how
much the spline $\sigma$ differs from a trajectory of the vector 
field. The spline satisfies $\sigma'(t) = F(\sigma(t))$ at the 
collocation points $t_{j-1}, t_j$ and $(t_{j} + t_{j-1})/2 $. 
If $\gamma$ is the trajectory of the vector field through one
these points, this implies that 
$\sigma - \gamma = O(|t_{j} - t_{j-1}|^4)$. Since this classical
estimate is based upon the assumption that the norm of the vector 
field is $O(1)$, it is only likely to hold for intervals that are
short on the fast time scale. 
Globally, the trajectories of the flow display a strong separation 
due to the normal hyperbolicity. In \emph{Fenichel
coordinates}~\cite{Jones}, stable coordinates converge rapidly to the slow
manifold while unstable coordinates diverge rapidly from the
slow manifold. In the case of a one dimensional slow manifold, 
\emph{shadowing}~\cite{Bowen} implies that any \emph{pseudo-trajectory}
pieced together from local approximations to the flow will lie close
to a unique trajectory of the flow. Moreover, in this case, different
choices of boundary conditions that lie in the same strong stable 
manifold at $a$ and the same strong unstable manifold at $b$ yield 
trajectories that are exponentially close to each other and to the 
slow manifold outside of small subintervals near the ends of the time
interval $[a,b]$. Consequently, the value of $F$ will be $O(1)$ on the
slow time scale and solutions of $E$ are expected to give 
quite accurate approximations to the slow manifold. Rather than 
pursuing more careful theoretical analysis of the algorithm, this 
paper tests its implementation on several examples.

\section{Examples}

This section presents four examples:
\begin{enumerate}
\item
A linear system for which there are explicit solutions of both the solutions of
the differential equations and the boundary value solver,
\item 
A three dimensional version of the Morris-Lecar model for bursting
neurons that was used by David Terman in his analysis of the
transition between bursts with different numbers of spikes~\cite{Ter,Ter1},
\item
A three dimensional system whose homoclinic orbits yield traveling
wave profiles for the Fitzhugh-Nagumo model~\cite{Champ},
\item
A four dimensional model of two coupled neurons studied by
Guckenheimer, Hoffman and Weckesser~\cite{GHW}.
\end{enumerate}

\subsection{Slow manifolds of a linear system}
The general solution of the linear vector field
\begin{equation}
\label{elinear}
\begin{split}
\eps \dot{x}_1 & =  y-x_1  \\
\eps \dot{x}_2 & =  x_2 \\
\dot{y} & = 1
\end{split}
\end{equation}
is 
$$(x_1,x_2,y)(t) = (y(0)-\eps+t+(x_1(0)-y(0)+\eps)\exp(-t/\eps),x_2(0)\exp(t/\eps),y(0)+t).$$
This explicit solution provides a benchmark for evaluating the accuracy of the algorithm
described above.
The slow manifold of the system is the line $y=x_1+\eps,x_2=0$ containing the trajectories
$(x_1,x_2,y)(t) = (y(0)-\eps+t,0,y(0)+t)$. 

The discretized equations of the algorithm 
can also be solved explicitly for system~\eqref{elinear}. The first step in doing so
is to observe that the equations for $x_2$ and $y$ are
separable from those for $x_1$, and this remains the case for the
discretized equations of the boundary value solver. Substituting the equations 
for the $y$ variable into the boundary value equations produces the equation 
$y_{j+1}-y_j = t_{j+1}-t_j$ on each mesh interval. If a boundary condition is 
imposed on one end of the time interval $[a,b]$, these equations yield a solution 
that is a discretization of an exact solution of the differential equation. 
Convergence occurs in a single step. 

Assume now that $y_{j+1}-y_j = t_{j+1}-t_j$ and set $w_j = y_j - (x_1)_j - \eps$ to be 
he difference between the $x_1$ coordinate of a point and a point of the slow manifold. 
The boundary value  equations become
$$ \frac{\delta^2-6\delta\eps+12\eps^2}{8\delta\eps^2}w_{j} - 
\frac{\delta^2+6\delta\eps+12\eps^2}{8\delta\eps^2}w_{j+1} = 0$$ 
with $\delta = y_{j+1}-y_j = t_{j+1}-t_j$. Note that these equations are satisfied 
when the $w_j$ vanish, so discretizations of exact solutions along the slow manifold 
satisfy the boundary value equations. Solving the equation for $w_{j+1}$ in terms 
of $w_j$ yields
$$ w_{j+1} = \frac{\delta^2-6\delta\eps+12\eps^2}{\delta^2+6\delta\eps+12\eps^2}w_{j}$$
Like the solutions of the differential equation, the values ${w_j}$
decrease exponentially as a function of time. The ratio $\rho_j = w_{j+1}/{w_j}$ is a 
function of $(\delta/\eps)$ whose Taylor expansion agrees with that of $\exp(-\delta/\eps)$ 
through terms of degree $4$, and its value always lies in the interval $(0,1)$. Thus 
the the solutions of the boundary value equation converge geometrically toward the 
slow manifold along its stable manifold with increasing time. If the mesh 
intervals have length $\delta \le \eps$, then the relative error of the decrease satisfies
$$0 < \frac{\rho_j(\frac{\delta}{\eps}) - \exp(\frac{\delta}{\eps})}{\exp(\frac{\delta}{\eps})} < 0.0015$$
For large values of $\delta/\eps$, the solution is no longer accurate near $t=a$ if 
the boundary conditions do not satisfy $y_0=(x_1)_0+\eps$. A similar, but simpler argument
establishes that the solution of the discretized problem converges to the slow manifold 
at an exponential rate with decreasing time from $t=b$. Thus, the boundary value solver 
is stable and yields solutions that qualitatively resemble the exact solution for all 
meshes when applied to this linear problem. In particular, the solution of the discretized
problem is exponentially close to the slow manifold away from the ends of the time interval
$[a,b]$. As the mesh size decreases to zero, the algorithm has fourth order 
convergence to the exact solution. 

\subsection{Bursting Neurons}

\begin{figure}[!htb]
  \begin{center}
\subfigure[]{
   \includegraphics[height=2.75in]{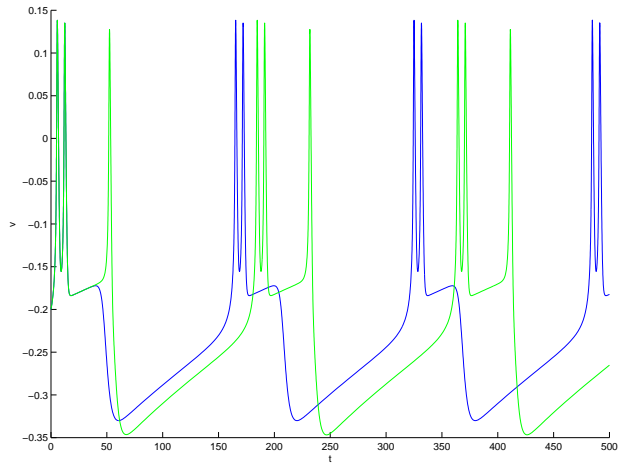}}
\subfigure[]{
    \includegraphics[height=2.75in]{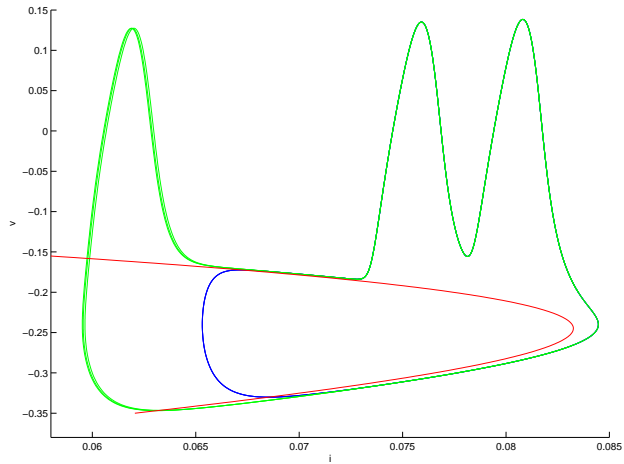}}
 \end{center}
  \caption{(a) Time series of $v$ showing two different periodic orbits of system~\eqref{mleqn}.
Parameter values are $(k,\eps) = (-0.24,0.00412234944)$. The blue orbit has two spikes per burst; 
the green orbit has three spikes per burst. (b) Phase portraits of the same periodic orbits 
projected onto the $(I,v)$ plane. A part of the critical manifold $S$ is shown in red.}
\label{terman_traj}
\end{figure}

Action potentials are a primary means for communicating information within the nervous 
system~\cite{KS}. Neurons are said to burst~\cite{Gorman} when they fire several consecutive 
action potentials between ``silent'' periods free of bursts. There is no universally 
accepted definition of bursts, but computational models are widely used to predict in 
terms of membrane channel properties when a neuron will burst. Rinzel~\cite{Rinzel} 
introduced a singular perturbation perspective to the investigation of bursting in model 
neurons, viewing the phenomenon as a relaxation oscillation in which a system makes fast 
time scale transitions between slowly varying equilibrium and periodic attractors. Several 
classifications of bursting distinguish qualitatively different dynamics. For example,
Izhikevich~\cite{Iz} classifies bursts in terms of the bifurcations that mark the 
transitions between bursts and silent intervals. Terman~\cite{Ter1} studied changes in 
the number of spikes per burst that occur as system parameters are varied. He gave numerical 
examples in a version of the Morris-Lecar model~\cite{ML} first analyzed by Rinzel and 
Ermentrout~\cite{RE}:
\begin{equation}
\label{mleqn}
\begin{split}
   \dot{v} & = I - 0.5(v+0.5)-2w(v+0.7)-0.5(1+\tanh(\frac{v-0.1}{0.145})(v-1)  \\
   \dot{w} & = 1.15(0.5(1+\tanh(\frac{v+0.1}{0.15})-w)\cosh(\frac{v-0.1}{0.29}) \\
   \dot{I} & = \eps(k-v)
\end{split}
\end{equation}
This system has periodic bursting solutions with different numbers of spikes per burst as 
the parameters $\eps$ and $k$ vary. Figure~\ref{terman_traj} illustrates that there are 
narrow parameter ranges with two stable periodic orbits having different spike numbers. 
Terman described the dynamics of the transition from periodic solutions with $n$ 
spikes to those with $n+1$, relying upon numerical simulations of trajectories in his
analysis. Flow along a slow manifold of saddle type is a central aspect of this transition,
but the trajectory simulations are incapable of following trajectories that remain close
to this slow manifold for more than a short distance. The boundary value solver
introduced in this paper is used to compute trajectories that contain 
segments which follow the slow manifold of saddle type. Geometric structures involved
with the transition from $n$ to $n+1$ spikes are visualized, 
and the analysis of the transition from $n$ to $n+1$ spikes is carried further. 

The fast subsystem of \eqref{mleqn} is the Morris-Lear model for action potentials 
of barnacle muscle~\cite{ML}. The Morris-Lecar model itself has a rich dynamical structure~\cite{RE}. 
There is an interval of values for $I$ in which the system has three equilibrium points. 
Saddle-node bifurcations occur at the endpoints of this interval. The equilibrium 
points of the Morris-Lecar model constitute the critical manifold of \eqref{mleqn}, 
and its saddle-node bifurcations (with varying $I$) are the folds of the critical 
manifold. There is also a family of periodic orbits that collapses at a subcritical 
Hopf bifurcation and terminates at a homoclinic bifurcation. This family of periodic 
orbits is folded, i.e., there is a saddle-node of limit cycle bifurcation~\cite{GH} within the 
family. The periodic orbits of the family between the fold and homoclinic bifurcations are stable. 

\begin{figure}[!htb]
  \begin{center}
\includegraphics[height=5in]{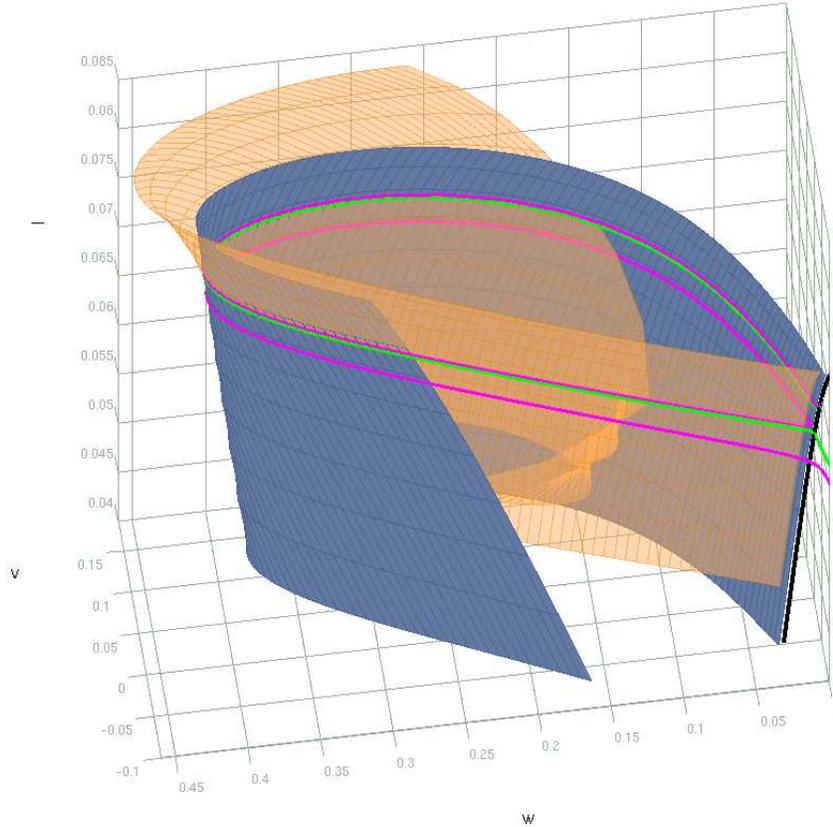}
 \end{center}
  \caption{Stable (blue) and unstable (orange) manifolds of the slow manifold (black) 
of saddle type in system~\eqref{mleqn} showing an intersection close to the homoclinic 
orbit of the singular limit of this system. The green and magenta curves are two 
trajectories with initial conditions that lie on opposite sides of the intersection. 
Parameters are $(k,\eps) = (-0.22,0.002)$.}
\label{terman_mflds}
\end{figure}

The bursting orbits of the vector field \eqref{mleqn} follow a branch 
of the critical manifold of equilibrium points to one of its folds, jump to the family 
of stable periodic orbits, follow this family to its homoclinic bifurcation and then 
jump back to the branch of stable equilibria. These bursting orbits occur when the 
value of the parameter $k$ is chosen so that $I$ increases slowly during the quiescent 
part of the cycle and decreases slowly during the active spiking portion of the cycle. 
See Figure~\ref{terman_traj}. The homoclinic orbit of the singular limit $\eps=0$ is a 
transversal intersection of the stable and unstable manifolds of the branch of saddle 
equilibria of ~\eqref{mleqn}. The branch of equilibria become a slow manifold $S$ of 
saddle type when $\eps>0$ and the homoclinic orbit persists as an intersection of the 
stable and unstable manifolds $W^s(S),W^u(S)$ of $S$. The transition between $n$ and 
$n+1$ spikes per burst occurs when the periodic bursting cycle encounters the 
intersection of $W^s(S)$ and $W^u(S)$. The final spike of a periodic orbit with $n+1$ 
spikes follows the intersection of $W^s(S)$ and of $W^u(S)$ back to $S$ before jumping 
to the stable slow manifold. Figure~\ref{terman_mflds} visualizes $S,W^s(S)$ and 
$W^u(S)$, and shows two trajectories that bracket the intersection of $W^s(S)$ and 
$W^u(S)$. 

Because the system \eqref{mleqn} is smooth and does not have an equilibrium point 
near the intersection of $W^s(S)$ and $W^u(S)$, the transition from $n$ to $n+1$ spikes 
consists of trajectories that undergo a continuous evolution. 
These trajectories contain saddle canards, 
segments that follow $S$ for varying lengths of time before leaving $S$ along one of its 
strong unstable manifolds. Trajectories lying close enough to $W^s(S)$ turn and flow along 
$S$ when they approach it. The distance that they travel along $S$ before leaving along 
its unstable manifold $W^u(S)$ depends logarithmically on the initial distance of the 
trajectory to $W^s(S)$. If close enough, the trajectory will follow $S$ all the way to its 
end near a fold of the critical manifold before making a fast excursion to the stable 
slow manifold. Accurate computation of $S$ is essential to understanding the details of the 
transition from bursts with $n$ spikes to bursts with $n+1$ spikes. 

\begin{figure}[!htb]
  \begin{center}
\includegraphics[height=5in]{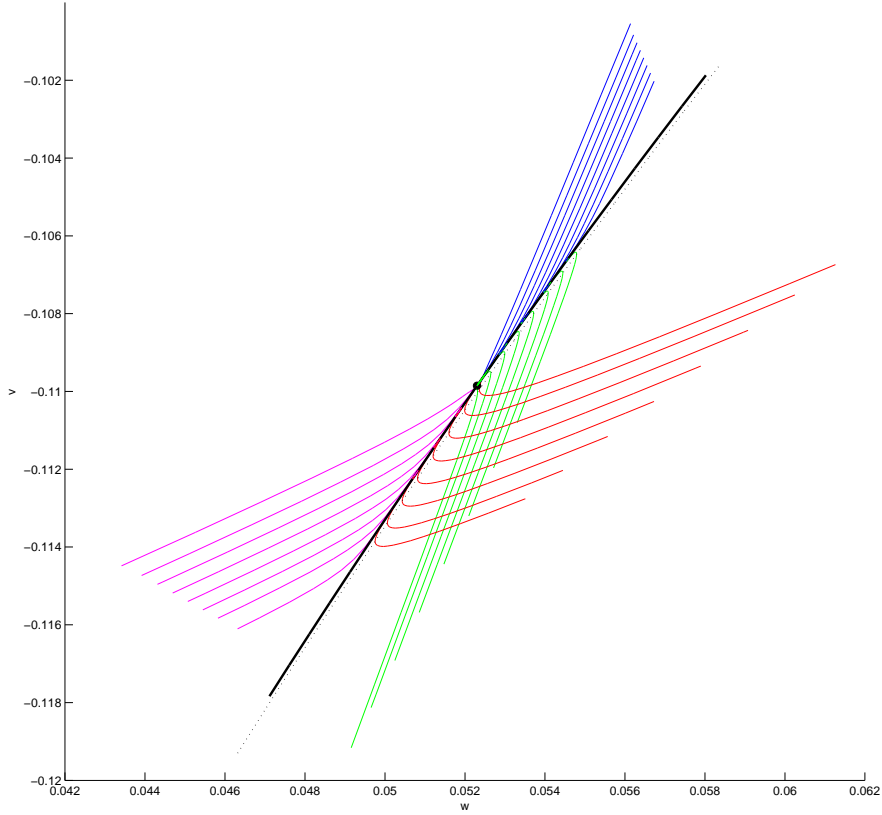}
 \end{center}
  \caption{Trajectories with initial conditions close to the slow manifold test the accuracy
of the slow manifold computations. The slow manifold is drawn as a heavy black curve and the 
critical manifold is drawn as a dotted black curve. Initial conditions for thirty-two trajectories
are chosen at distances $10^{-4},10^{-5},10^{-6},10^{-7},10^{-8},10^{-9},10^{-10},10^{-11}$
along the strong stable and unstable manifolds at the point 
$(-0.109854033586602,0.052299738361417,0.025187193494031)$ on the slow manifold, which is 
drawn as a filled black circle. The trajectories are computed for a time interval $\pm 0.01$ and color
coded so that the trajectories along the two branches of the strong unstable manifold are 
drawn blue and green while the trajectories along the strong stable manifold are drawn red and 
magenta. Parameter are $(k,\eps) = (-0.22,0.002)$ and the objects are projected into the $(w,v)$ 
plane.}
\label{terman_test}
\end{figure}

The critical manifold of the vector field~\eqref{mleqn} is given by explicit formulas 
when  parametrized by $v$. Uniform meshes of values for $v$ were used to generate starting
values for the boundary value computation of the slow manifold $S$. To compute $S$,
the vector field~\eqref{mleqn} was rescaled so that $\dot{I} = 1$. With this rescaling,
the value of $I$ remains constant during the Newton iteration to find the solutions. 
Typical meshes that do not come close to the fold points result in convergence of Newton's 
method within three or four steps. Figure~\ref{terman_test} illustrates the accuracy of 
the computations of $S$ and the behavior of numerical simulations of trajectories that start
near $S$. A point $p$ on $S$ is chosen, and the Jacobian of the fast subsystem at this point is 
computed to obtain approximations for the directions of its strong stable and unstable 
manifolds. If $p$ does lie on the slow manifold, then trajectories with initial conditions 
on opposite sides of $S$ on its strong unstable manifold will flow along $S$ but then jump 
in opposite directions. Similarly, backward trajectories with initial conditions 
on opposite sides of $S$ on its strong stable manifold will flow along $S$ but then jump 
in opposite directions. If $p$ is displaced from $S$, its distance to $S$ can be estimated
by finding the closest pairs of bracketing trajectories that do jump from $S$ in opposite 
directions. Figure~\ref{terman_test} displays the results of such a test. Eight pairs of 
trajectories displaced along the strong unstable manifold at distances $10^{-k},4 \le k \le 11$
are plotted in blue and green, and eight pairs of backward
trajectories displaced along the strong stable manifold at distances $10^{-k},4 \le k \le 11$
are plotted in red and magenta. Pairs of trajectories displaced by distance $10^{-12}$ (not drawn
in Figure~\ref{terman_test}) fail the test, jumping in the same direction. This suggests that 
the distance from $p$ to the slow manifold
is smaller than $10^{-11}$. Note also that increments in the distance that each successive pair of 
bracketing trajectories flows along $S$ are similar, consistent with the exponential separation of 
trajectories within the strong stable and unstable manifolds. Extrapolating these increments yields 
the estimate that a numerically simulated trajectories starting on the slow manifold near $p$ will
only be able to remain close to $S$ for time approximately $0.01$. This estimate is based on 
round-off error of the order of $10^{-16}$ and the observation that the times at which trajectories
displaced from $p$ by distances $10^{-9}$ and $10^{-11}$ appear to jump from $S$ are approximately
$0.003$ and $0.005$. These crude estimates explain why initial value solvers are unable to follow
the continuous evolution of trajectories in the transition from $n$ to $n+1$ spikes per burst.
The value of $v$ at $p$ is approximately $-0.11$ and the jump from $S$ of numerically simulated
trajectories seems to occur before $v$ increases to $-0.1$, but the fold of the critical manifold 
occurs when $v$ is approximately $-0.034$. The exponential instability of $S$ in both forward and
backward directions precludes initial value solvers from computing trajectories that flow along
$S$ from the intersection of $W^s(S)$ and $W^u(S)$ to the fold of $S$.

Computation of periodic orbits with long canard segments near the slow manifold of saddle type 
appears to be challenging, even with continuation methods. Computation of the slow manifold $S$ 
with the boundary value solver introduced here can be coupled with the analysis of Terman~\cite{Ter} 
and Lee and Terman~\cite{Ter1} to solve this problem. Computations of the slow manifold $S$ are 
augmented with numerical forward and backward simulations of trajectories that terminate at a 
cross-section along the family of periodic orbits. This extends the approach introduced by 
Guckenheimer and Lamar~\cite{GL} to efficiently compute periodic orbits containing canards. 

\begin{figure}[!htb]
  \begin{center}
\includegraphics[height=5in]{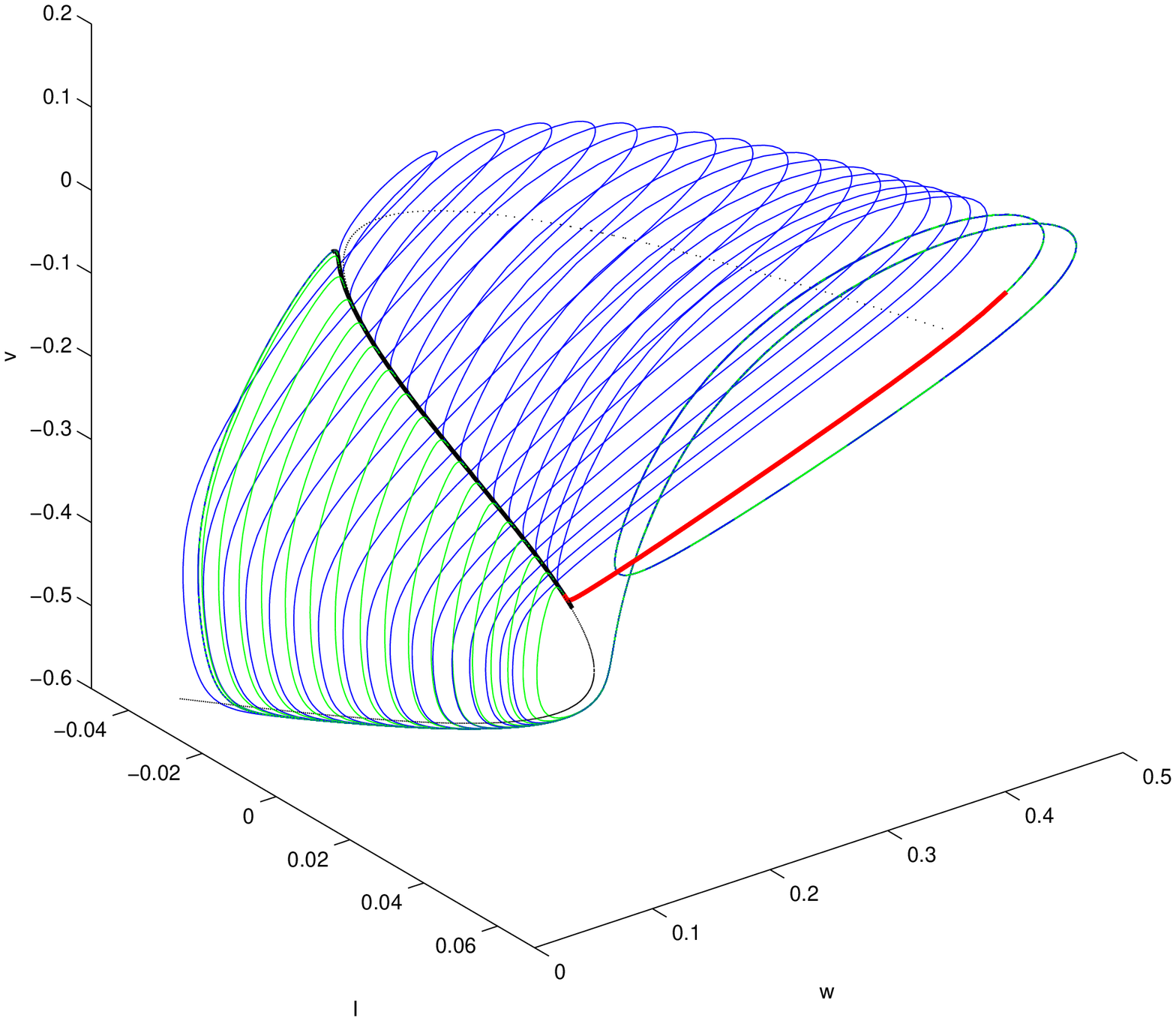}
 \end{center}
  \caption{The unstable manifold of the slow manifold $S$ of saddle type in system~\eqref{mleqn} is swept out by trajectories starting close to the slow manifold. The blue and green trajectories are followed until they intersect the cross-section $I=0.075$ with $I$ increasing. Red trajectories on the stable manifold of $S$ are followed backward until they intersect this cross-section. The trajectories that are drawn reach the cross-section close to the points on the unstable manifold of $S$. Parameters are $(k,\eps) = (-0.22,0.006366)$.}
\label{terman_umfld}
\end{figure}

\begin{figure}[!htb]
  \begin{center}
\subfigure[]{
   \includegraphics[height=1.75in]{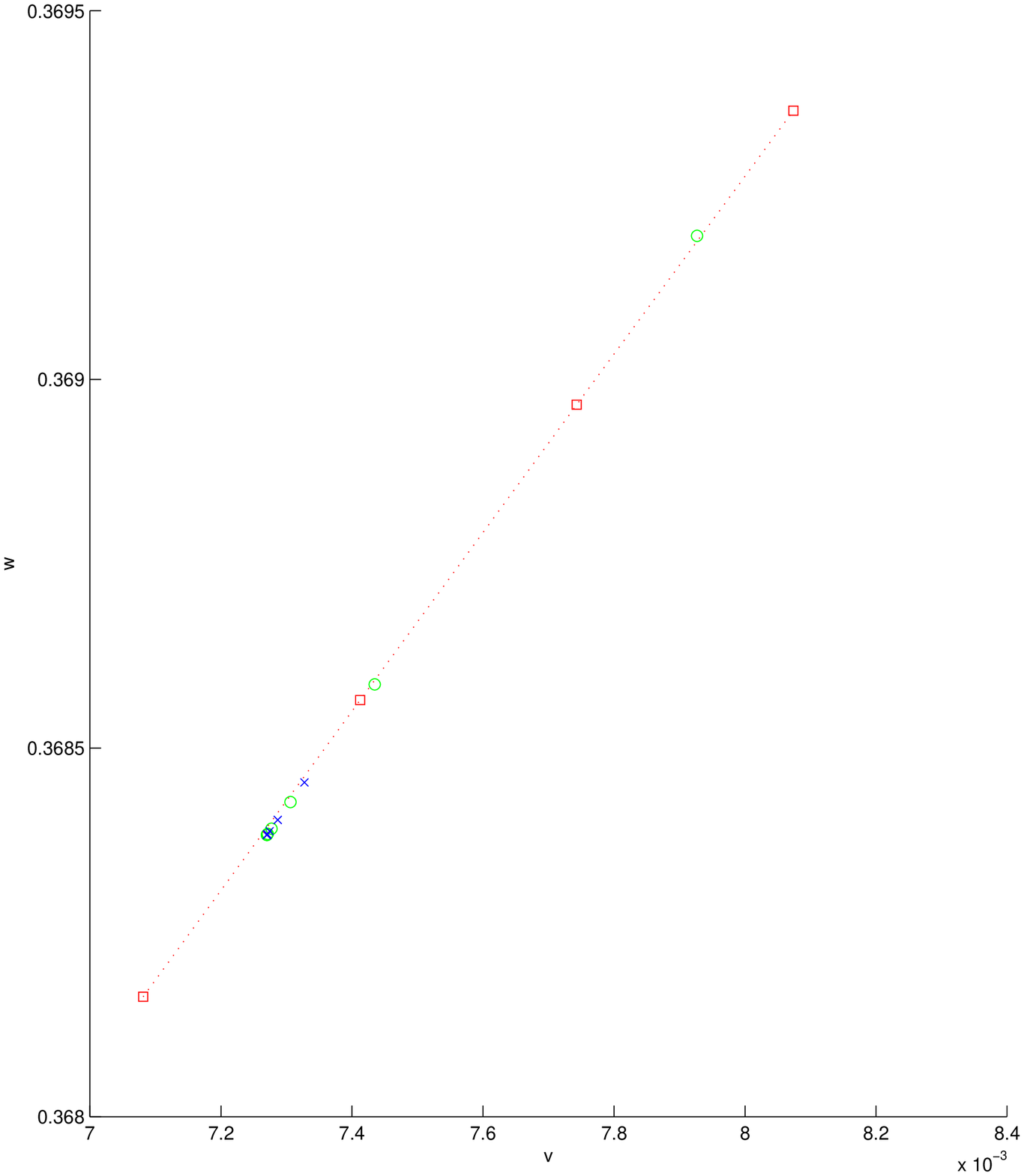}}
\subfigure[]{
   \includegraphics[height=1.75in]{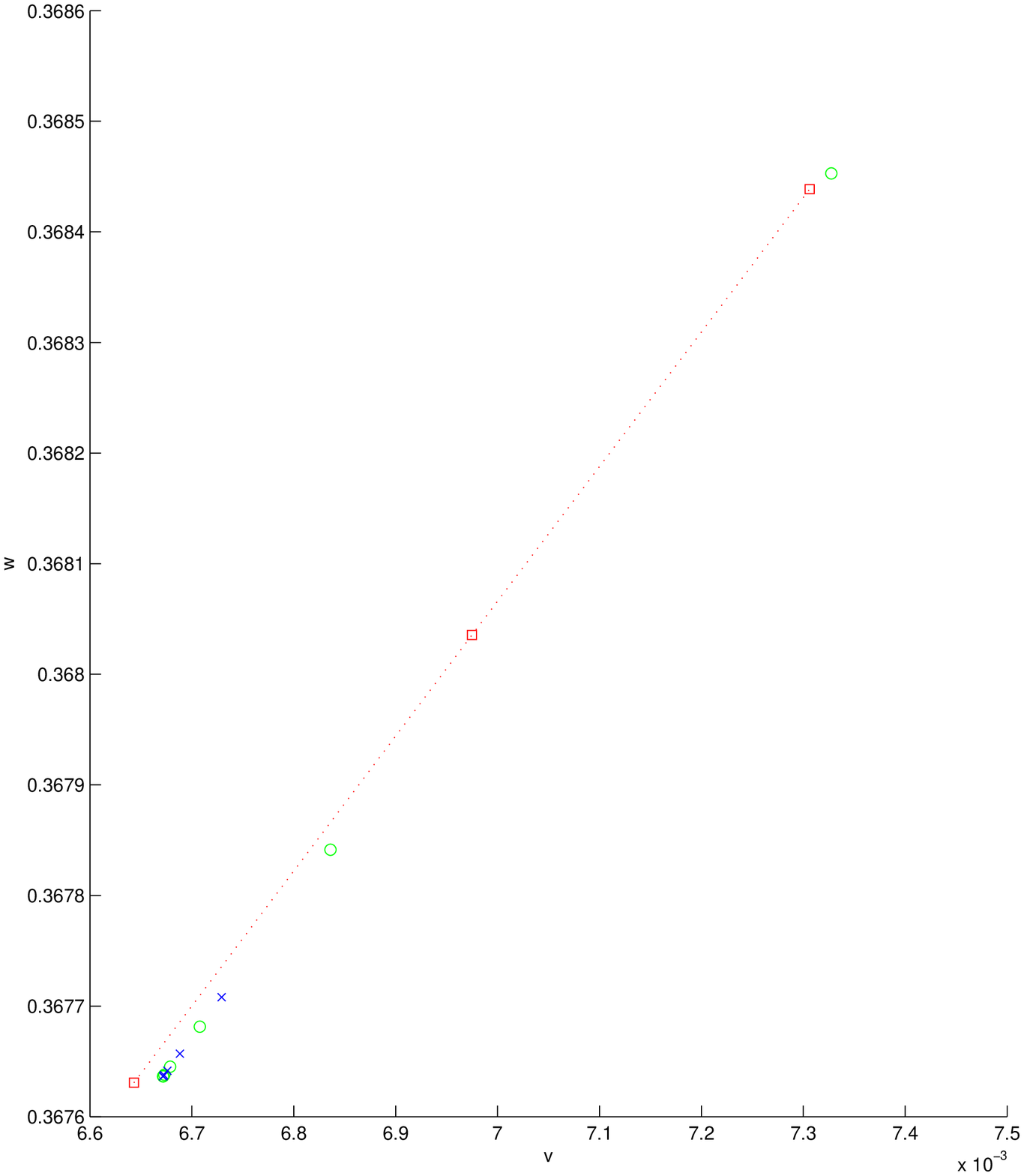}}
\subfigure[]{
   \includegraphics[height=1.75in]{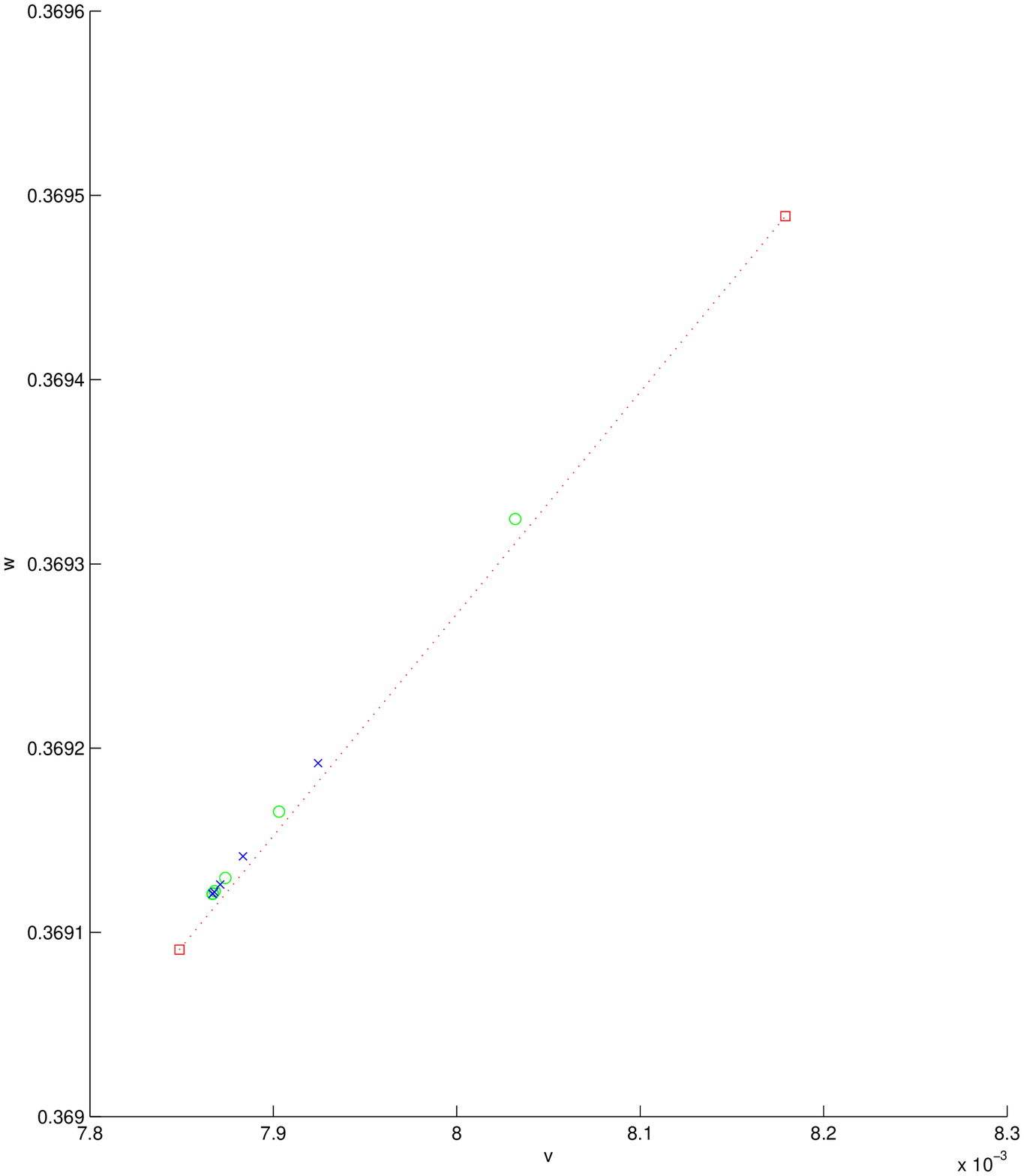}}\
 \end{center}
  \caption{(a) Endpoints of the trajectories displayed in Figure~\ref{terman_umfld} with the plane $I=0.075$. Blue crosses and green circles are used to mark points on $W^u(S)$; red squares denote points on $W^s(S)$. The red dotted curve is a piecewise linear connection between the points on $W^s(S)$, showing that $W^u(S)$ and $W^s(S)$ almost intersect for these parameter values $(k,\eps) = (-0.22,0.006366)$. (b) A similar plot to (a) for parameter values $(k,\eps) = (-0.22,0.006362)$. The points of $W^u(S)$ lie below those of $W^s(S)$. (c) A similar plot to (a) for parameter values $(k,\eps) = (-0.22,0.006367)$. The points of $W^u(S)$ lie above those of $W^s(S)$.}
\label{terman_sect}
\end{figure}

Figure~\ref{terman_umfld} visualizes the invariant manifold $W^u(S)$ as a collection of trajectories for parameter values $(k,\eps) = (-0.22,0.006366)$ in system~\eqref{mleqn}. The heavy black curve is a segment of the slow manifold $S$ of saddle type, and the dotted black curve is the critical manifold. At twenty initial points along $S$, trajectories have been computed with initial conditions displaced from $S$ along its strong unstable manifolds by a distance $0.00005$. The trajectories starting on one side of $S$ are drawn blue and the trajectories starting on the other side of $S$ are drawn green. The blue  trajectories make a loop around the unstable branch of the slow manifold and then flow past $S$ to the stable branch of the slow manifold. The green branches flow to the stable branch of the slow manifold with $v$ decreasing. Both sets of branches then turn and flow along the stable branch of the slow manifold. When they reach the fold of the slow manifold, they jump to the family of rapid oscillations. As trajectories follow these oscillations, $I$ decreases. The displayed trajectories are terminated when they reach the plane $I = 0.075$ with $I$ decreasing. The red curves displayed in Figure~\ref{terman_umfld} are four backwards trajectories that begin at distance $5 \times 10^{-8}$ from $S$ along its stable manifold and end on the cross-section $I = 0.075$. These trajectories were chosen on a short section of $S$ so that they reach the cross-section $I = 0.075$ near the ends of the blue and green trajectories. Figure~\ref{terman_sect}(a) shows the ends of the blue, green and red trajectories with the cross-section $I = 0.075$. Figure~\ref{terman_sect}(b) and (c) show similar plots for the system with parameter values $(k,\eps) = (-0.22,0.006362)$ and $(k,\eps) = (-0.22,0.00637)$. As $\eps$ varies, these plots demonstrate that the trajectories in the unstable manifold of $S$ sweep across the stable manifold of $S$. 

\begin{figure}[!htb]
  \begin{center}
\subfigure[]{
   \includegraphics[height=2.75in]{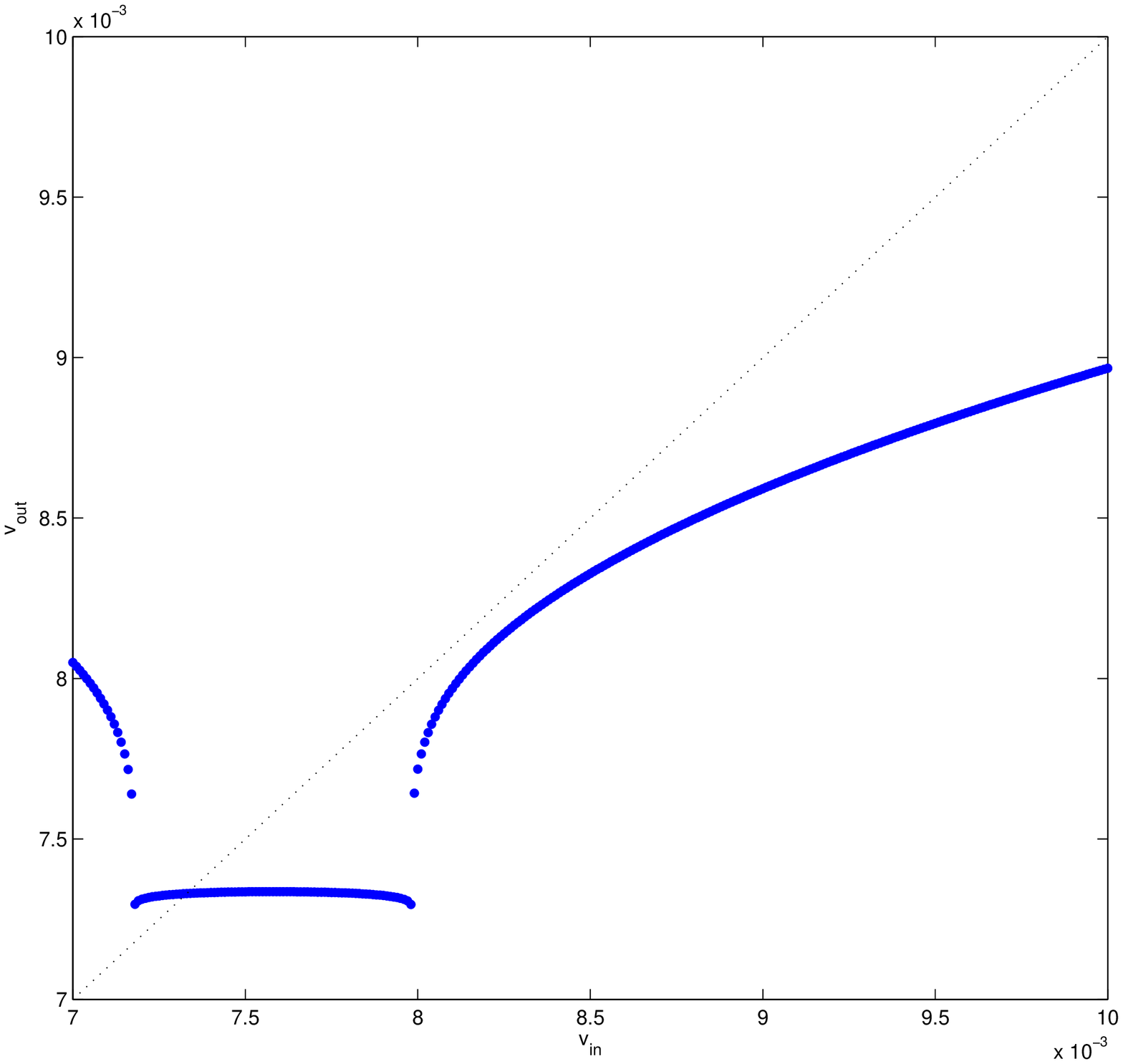}}
\subfigure[]{
   \includegraphics[height=2.75in]{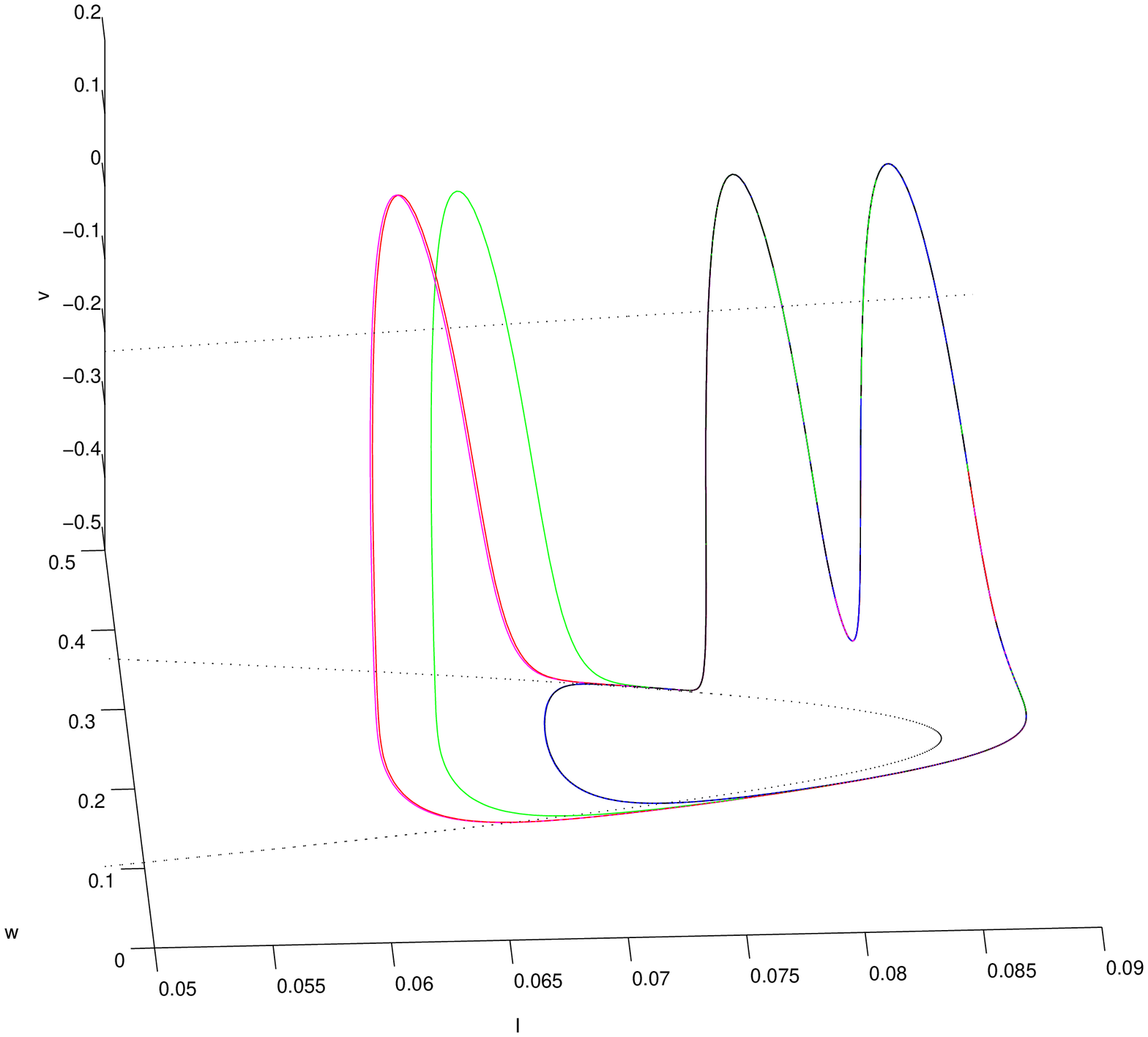}}
 \end{center}
  \caption{(a) Return map of the segment along the line $w = 1.2107v+0.35959, v \in [0.007,0.01], I = 0.075$. Three hundred initial points were chosen along this segment and computed until they return to $I = 0.075$. The axes are initial and final values of $v$. Parameter values are $(k,\eps) = (-0.22,0.006366)$. (b) Five trajectories among the three hundred computed for the return map in (a). Four of these trajectories bracket the jumps of the return map; the fifth has initial condition at the local maximum of the return map.}
\label{terman_retmap}
\end{figure}

Figure~\ref{terman_umfld} supports the following procedure for finding periodic orbits containing canards. Fix a short segment $\Sigma$ transverse to $W^u(S)$. With varying $\eps$, trajectories with initial conditions on $\Sigma$ sweep out a three dimensional manifold $M$ in $(v,w,I,\eps)$ space. The Exchange Lemma~\cite{JK} implies that if $M$ intersects $W^s(S)$ transversally in $(v,w,I,\eps)$ space, then part of $M$ will stretch along the length of $S$ and depart from it along $W^u(S)$. In particular, $M$ will intersect $\Sigma$, giving a unique value of $\eps$ for which there is a periodic orbit intersecting $\Sigma$. Figure~\ref{terman_sect} gives numerical evidence that $M$ does intersect $W^u(S)$, and it indicates that the value of $\eps$ will be almost constant along the family of periodic orbits containing canards. Computing trajectories with initial conditions on $\Sigma$ with an initial value solver will not produce these periodic orbits. Figure~\ref{terman_retmap}(a) shows a return map, giving initial and final values for the variable $v$, with 300 initial points chosen on a linear approximation to the intersection of $W^s(S)$ with $I=0.075$ shown in Figure~\ref{terman_sect}(a). This return map has two apparent jumps. The trajectories beginning between the two jumps make three spikes before returning to $I=0.075$ while the other trajectories make two spikes before returning. None of the trajectories flows along $S$ with $I$ decreasing to a value smaller than $0.065$. Figure~\ref{terman_retmap}(b) plots five of the trajectories from the return map, four that bracket the jumps and one from the local maximum of the return map in Figure~\ref{terman_retmap}(a).

The approximations of $S$ obtained with the boundary value solver can be used in 
approximating periodic orbits with canards. The strategy we propose is illustrated 
by Figures~\ref{terman_umfld} and \ref{terman_sect}. The periodic orbit will be calculated in three segments that are illustrated as black, red and blue/green curves in Figure~\ref{terman_umfld}. The canard segment of a periodic orbit is exponentially close to $S$ except at its arrival and departure points. When the periodic orbit departs from $S$, it will follow a blue or green trajectory starting at a point exponentially close to $W^u(S)$ that is numerically indistinguishable from points on $W^u(S)$. Similarly, the orbit segment that arrives at $S$, does so at a point that is exponentially close to $W^s(S)$, that is numerically approximated by backward integration beginning at a point of $W^s(S)$. As $\eps$ (or another parameter) is varied, the forward trajectory along $W^u(S)$ and backward trajectories along 
 $W^s(S)$ sweep out a curve and a surface of intersection with a cross-section in $(v,w,I,\eps)$ space (here $I=0.075$). A root solver can be used to locate a parameter value for which a trajectory of $W^u(S)$ and one on $W^s(S)$ arrive at the same point of the cross-section. The periodic orbit will then be approximated by the union of the two trajectories and a curve that flows along $S$ from the chosen arrival point to the chosen departure point. Normal hyperbolicity implies that there is a unique trajectory that connects these two points. As shown in Figure~\ref{terman_sect}, the intersections occur for $\eps \approx 0.006366$ for all arrival and departure points.

\begin{figure}[!htb]
  \begin{center}
\includegraphics[height=3in]{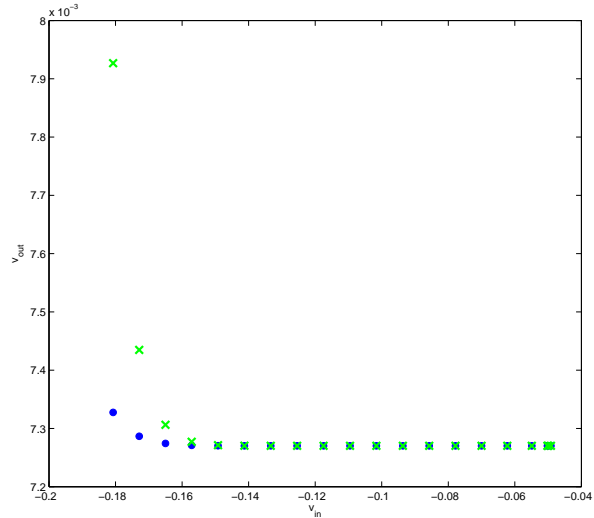}
 \end{center}
  \caption{Initial and final values of $v$ for trajectories in $W^u(S)$ 
ending in the plane $I=0.075$. Parameters are $(k,\eps) = (-0.22,0.006366)$.}
\label{terman_uv}
\end{figure}

Analysis of the bifurcations and attracting limit sets of the vector
field~\eqref{mleqn} requires additional information. Numerically, it
is necessary to ``fill in'' the jumps in the return map shown in
Figure~\ref{terman_retmap}, describing more carefully how the
trajectories with canards return and determining the stability of
trajectories containing canard segments.  Figure~\ref{terman_uv} plots
the final values of $v$ in trajectories on $W^u(S)$ vs. their initial
values of $v$. It is apparent that a large portion of $W^u(S)$
contracts enough when it flows along the stable branch of the slow
manifold that its intersection with the cross-section $I=0.075$ is
very small. The image appears to lie inside a disk centered at $(v,w)
= (0.0072701057,0.3683819196)$ of radius $10^{-10}$. The minimum
return value of $v$ in the points plotted in
Figure~\ref{terman_retmap} is approximately $0.007296$. Thus it
appears that the local minima of the return map are only a distance
about $3 \times 10^{-5}$ below the lowest points plotted in this
figure.

The variational equations of system~\eqref{mleqn} can be used to estimate how much expansion takes place along canard segments of trajectories and how much contraction takes place along the stable branch of the slow manifold. On points of the critical manifold with the same value of $I$, the strong unstable eigenvalue on the middle branch has larger magnitude than the weaker stable eigenvalue on the stable branch. If a canard segment is long enough, then the accumulated expansion will dominate the subsequent contraction that takes place on the stable branch of the slow manifold. This suggests that the return map of the system will have an expanding direction for canards that are sufficiently long, consistent with Figure~\ref{terman_retmap}. As the maximal canards of the return map move across the diagonal with changing parameters, the return map is likely to have chaotic invariant sets similar to those found in the Henon map~\cite{Henon}. The 
numerical computations reported here are insufficient to adequately determine the details of these invariant sets because the slow manifold is not computed close enough to the fold that that the maximal canards are determined with good accuracy.

Terman~\cite{Ter} suggests that the system~\eqref{mleqn} has trajectories with three different spike numbers in its bursts. The calculations here cast doubt abut whether this is possible
for trajectories that lie in the forward limit set of the system. 
For the parameters $(k,\eps) = (-0.22,0.006366)$, the decrease in the value of $I$
between spikes along the surface of oscillations appears to be approximately $0.0065$
in the region between the arrival of trajectories jumping from the fold of the stable branch
of the slow manifold to the intersection of $W^u(S)$ and $W^s(S)$. On the other hand,
the trajectories that flow along the stable branch of the slow manifold appear to pass 
by the fold in a set that has a diameter at least an order of magnitude smaller than the observed separation between spikes. Our analysis of canards makes it clear that the trajectories with long canard segments all flow through a tiny region as they pass the fold. lee and Terman~\cite{Ter1} give asymptotic estimates of the size of these regions in terms of $~\eps$ that 
also suggest that it is unlikely that the limit set of this system reaches the surface of oscillations in a set that is large enough to contain trajectories with three different spike numbers.

\subsection{Travelling Waves of the Fitzhugh-Nagumo Model}
\label{sec:fhn}

The FitzHugh-Nagumo equation is a model for the electric potential $u=u(x,\tau)$ of a nerve axon interacting with an auxillary variable $v=v(x,\tau)$ (see \cite{FitzHugh},\cite{Nagumo}):
\begin{equation}
\label{eq:fhn_original}
\left\{
\begin{array}{l}
\frac{\partial u}{\partial \tau }=\delta \frac{\partial^2 u}{\partial x^2}+f_a(u)-w+p \\
\frac{\partial w}{\partial\tau}  =\epsilon(u-\gamma w)
\end{array}
\right.
\end{equation}
where $f_a(u)=u(u-a)(1-u)$ and $p,\gamma,\delta$ and $a$ are parameters. Assuming a travelling wave solution with $t=x+s\tau$ to (\ref{eq:fhn_original}) we get:
\begin{eqnarray}
\label{eq:fhn_temp}
u'&=&v \nonumber\\
v'&=&\frac1\delta(sv-f_a(u)+w-p)\\
w'&=&\frac{\epsilon}{s}(u-\gamma w)\nonumber
\end{eqnarray}
A homoclinic orbit of (\ref{eq:fhn_temp}) corresponds to a travelling pulse solution in (\ref{eq:fhn_original}). An analysis of (\ref{eq:fhn_temp}) using numerical continuation has been carried out by Champneys et al. \cite{Champ}. They fixed the parameters $a=\frac{1}{10}$, $\delta=5$, $\gamma=1$ and investigated bifurcations in $(p,s)$-parameter space. We shall fix the same values and hence write $f_{1/10}(u)=:f(u)$. To bring (\ref{eq:fhn_temp}) into the standard form (\ref{sfs}) set $x_1:=u$, $x_2:=v$, $y:=w$ and change to the slow time scale:
\begin{eqnarray}
\label{eq:fhn}
\epsilon\dot{x}_1&=& x_2\nonumber\\
\epsilon\dot{x}_2&=&\frac15 (sx_2-x_1(x_1-1)(\frac{1}{10}-x_1)+y-p)=\frac15 (sx_2-f(x_1)+y-p)\\
\dot{y}&=&\frac{1}{s} (x_1-y) \nonumber 
\end{eqnarray}
We refer to (\ref{eq:fhn}) as ``the'' FitzHugh-Nagumo equation. Our goal is to use the fast slow structure of (\ref{eq:fhn}) and the SMST algorithm to compute its homoclinic orbits. The critical manifold $S$ of the FitzHugh-Nagumo equation is the cubic curve:
\begin{equation}
S=\{(x_1,x_2,y)\in\mathbb{R}^3: \quad x_2=0, \quad y=f(x_1)+p=:c(x_1)\}
\end{equation} 
The two local non-degenerate extrema of $c(x_1)$ yield the fold points of $S$. Denote the local minimum by $x_{1,-}$ and the local maximum by $x_{1,+}$. The critical manifold $S$ has three normally hyperbolic components:
\begin{equation*}
S_l=\{x_1<x_{1,-}\}\cap S, \quad S_m=\{x_{1,-}< x_1< x_{1,+}\}\cap S,\quad S_r=\{x_{1,+}<x_1\}\cap S
\end{equation*} 
Fenichel's theorem provides associated slow manifolds $S_{l,\epsilon}$, $S_{m,\epsilon}$ and $S_{r,\epsilon}$ outside neighbourhoods of the fold points. The manifolds $S_{l,\epsilon}$ and $S_{r,\epsilon}$ are of saddle-type for $\epsilon$ sufficiently small. The middle branch $S_{m,\epsilon}$ is completely unstable in the fast directions. Denote the unique equilibrium point of (\ref{eq:fhn}) by $q=(x_1^*,0,x_1^*)$. The location of $q$ depends on the parameter $p$ and $q$ moves along the cubic $S$. For the analysis of homoclinic orbits we shall assume that $q\in S_{l,0}$. In this case, the unstable manifold $W^u(q)$ is one-dimensional and the stable manifold $W^s(q)$ is two-dimensional. This also covers the case $q\in S_r$ by a symmetry in the FitzHugh-Nagumo equation and avoids the region where $q$ is completely unstable \cite{Champ},\cite{GuckCK1}. Homoclinic orbits exist if $W^u(q)\subset W^s(q)$.\\

\begin{figure}[htbp]
\centering
  \includegraphics[width=0.6\textwidth]{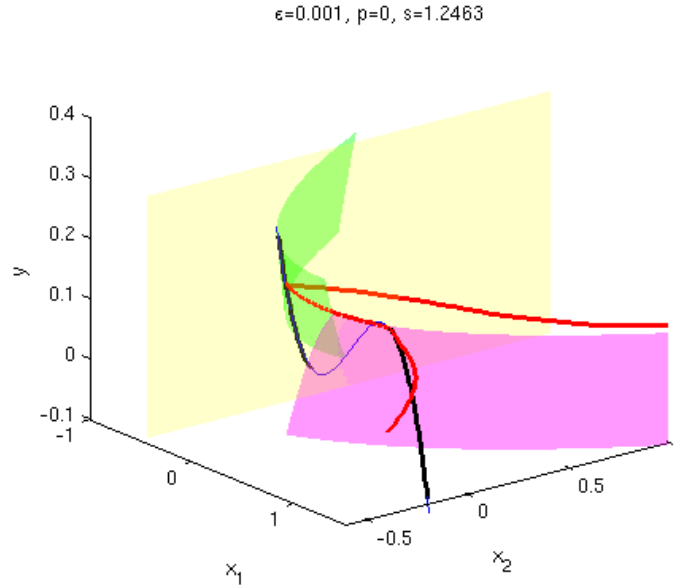} 
 \caption{Illustration of transversal intersection of stable and unstable manifolds of the slow manifolds $W^s(S_{l,\epsilon})$ (green) and $W^u(S_{r,\epsilon})$ (magenta). The manifolds are trucated at the yellow section $\Sigma$ and the trajectory $\gamma_{fw}\cup \gamma_{bw}$ started on $\Sigma$ at the transversal intersection point $x_{su}$ is shown in red.}
 \label{fig:hom_intersect_truncated}
 \end{figure}

We focus first on the case of relatively large wave speeds $s$ (``fast waves''). The existence proof of these homoclinic orbits contructs them as perturbations of a singular trajectory consisting of four segments: a fast subsystem heteroclinic connection from $q$ to $C_r$ at $y=x_1^*$, a slow segment on $C_r$, a fast subsystem heteroclinic from $C_r$ to $C_l$ at $y=x_1^*+c$ for some constant $c=c(p,s)>0$ and a slow segment on $C_l$ connecting back to $q$ \cite{JonesKopellLanger}. We aim to compute homoclinic orbits by a similar procedure for a given small $\epsilon>0$ in several steps:

\begin{enumerate}
 \item Find parameter values $(p_0,s_0)$ such that a homoclinic orbit exists very close or exactly at $(p_0,s_0)$. This can be achieved by a splitting algorithm without computing the homoclinic orbit, even for very small values of $\epsilon$ \cite{GuckCK1}. Carry out all the following compuations for $(p,s)=(p_0,s_0)$.
 \item \label{step1} Compute the slow manifolds $S_{\epsilon,l}$ and $S_{\epsilon,r}$ using the SMST algorithm.
 \item \label{step2} Compute the unstable manifold of the equilibrium $W^u(q)$ by forward integration.
 \item Define a section $\Sigma=\{x_1=c\}$ where the constant $c$ is chosen between $x_{1,-}$ and $x_{1,+}$ e.g. $c=(x_{1,-}+x_{1,+})/2$. Compute the transversal intersection of $W^s(S_{l,\epsilon})$ and $W^u(S_{r,\epsilon})$ on $\Sigma$, call the intersection point $x_{su}=(c,x_{2,su},y_{su})$ (see Figure \ref{fig:hom_intersect_truncated}). Integrate forward and backward starting at $x_{su}$ to obtain trajectories $\gamma_{fw}$ and $\gamma_{bw}$.
 \item \label{step4} To compute the homoclinic orbit we use the objects computed so far as approximants in different regions. Compute the closest points in $W^u(q)$ and $\gamma_{bw}$ to $S_{r,\epsilon}$ and concatenate $W^u(q)$ and $\gamma_{bw}$ to $S_{r,\epsilon}$ at these points. Proceed similarly with $S_{l,\epsilon}$ and $\gamma_{fw}$. Remove all parts of the slow manifolds not lying between the concatenation points and past the equilibrium $q$.   
\end{enumerate}

\begin{figure}[htbp]
 \subfigure[Slow manifolds $S_{l,\epsilon}$ and $S_{r,\epsilon}$ are shown in black and the unstable manifold of the equilibrium $W^u(q)$ is displayed in red.]{\includegraphics[width=0.5\textwidth]{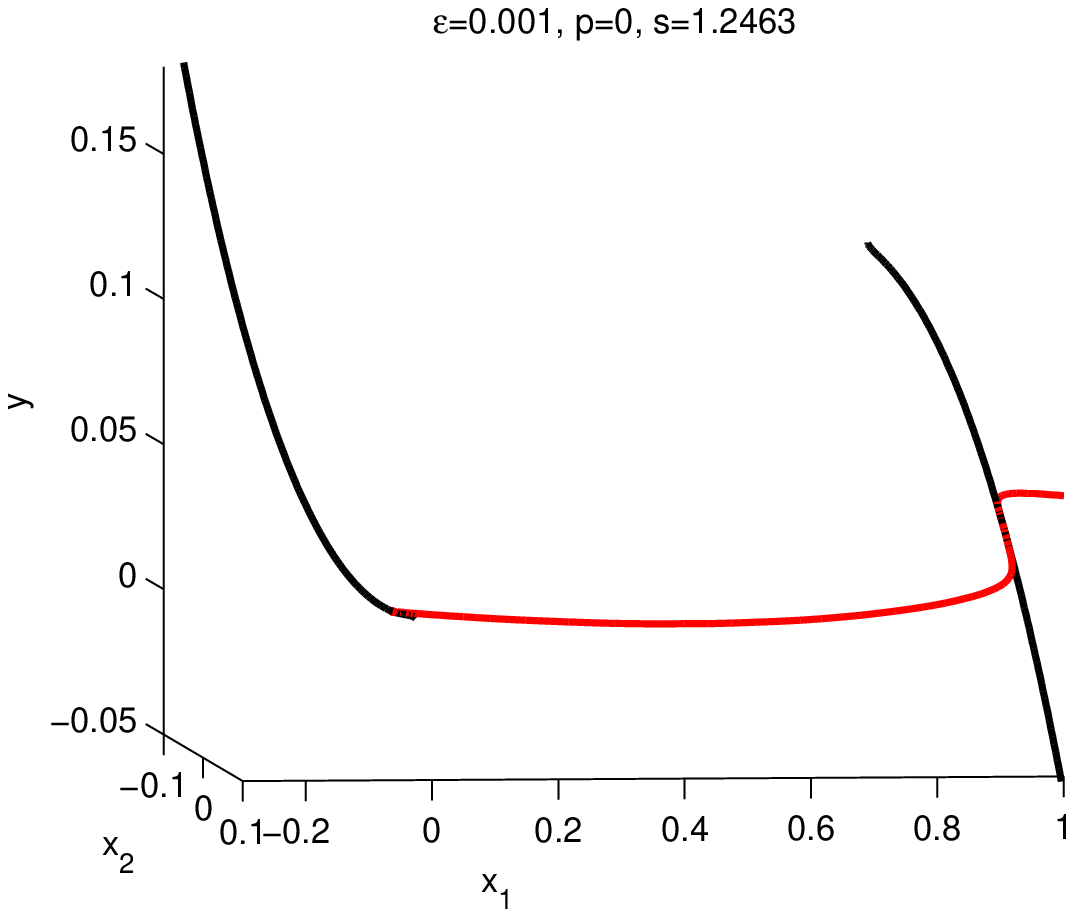}\label{fig:hom_umfld_smfld}}
 \subfigure[Pieces of the homoclinic orbit; slow segments in black, fast segments in red and $S$ shown in blue.]{\includegraphics[width=0.5\textwidth]{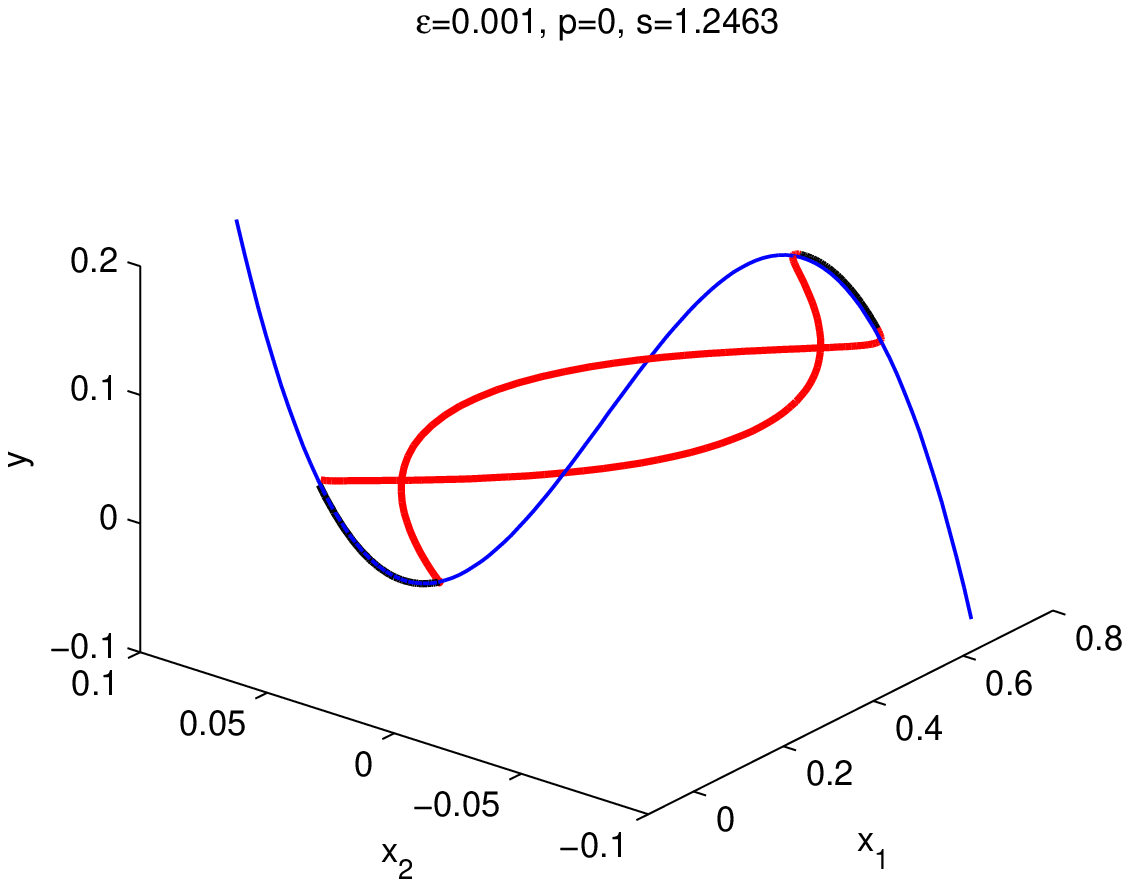}\label{fig:hom_orbit_fast_pieces}}
\caption{Illustration of the algorithm for computing homoclinic orbits in the FitzHugh-Nagumo equation.}
\label{fig:fast_red}
\end{figure}

\begin{figure}[htbp]
\centering
  \includegraphics[width=0.75\textwidth]{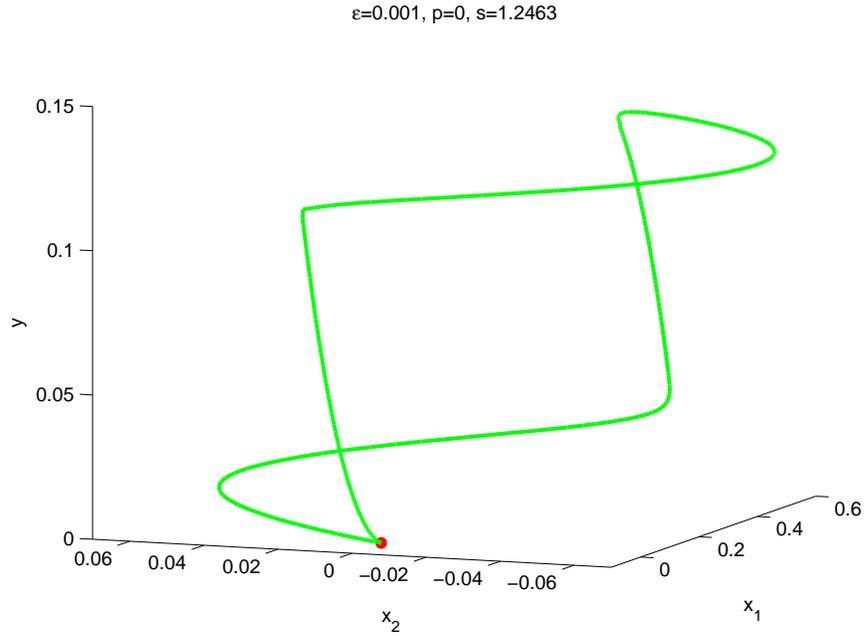} 
 \caption{Homoclinic orbit (green) of the FitzHugh-Nagumo equation representing a fast wave. The equilibrium point $q$ is shown in red.}
 \label{fig:hom_orbit_fast}
\end{figure}

Note that all figures for the fast wave case have been computed for $\epsilon=10^{-3}$, $p_0=0$ and $s_0\approx 1.2463$. This is a classical case \cite{JonesKopellLanger} for which the existence of homoclinic orbits is known. In Figure \ref{fig:hom_umfld_smfld} we show the result from the SMST algorithm and the unstable manifold of the equilibrium $W^u(q)$, i.e. the output of steps \ref{step1} and \ref{step2}. Due to the exponential separation along $S_{r,\epsilon}$ the trajectory $W^u(q)$ obtained from numerical integration cannot track the slow manifold for an O(1) distance and escapes after following the slow manifold for a very short time. This happens despite the fact that we have computed parameter values $(p_0,s_0)$ with maximal accuracy in double precision arithmetic at which we expect $W^u(q)$ to follow $S_{r,\epsilon}$ almost up to the fold point $x_{1,+}$. This observation is relevant to Figure \ref{fig:hom_orbit_fast_pieces} where the result of step \ref{step4} is shown. All the fast segments (red) had to be truncated almost immediately after they entered a neighourhood of a slow manifold. The final output of the algorithm after interpolation near the truncation points is shown in Figure \ref{fig:hom_orbit_fast}.\\  

\begin{figure}[htbp]
 \subfigure[``Single pulse'' homoclinic orbit]{\includegraphics[width=0.5\textwidth]{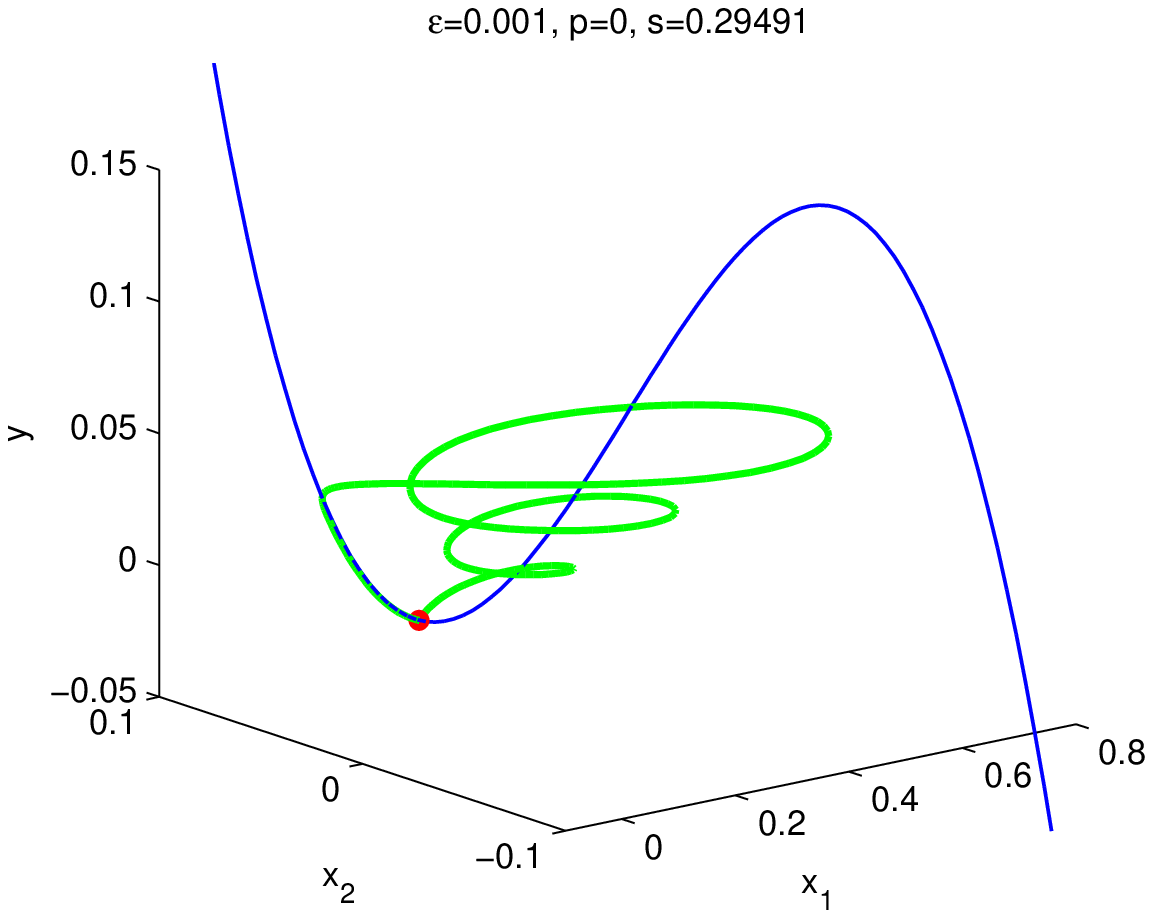}\label{fig:slow_wave1}}
 \subfigure[``Double pulse'' homoclinic orbit]{\includegraphics[width=0.5\textwidth]{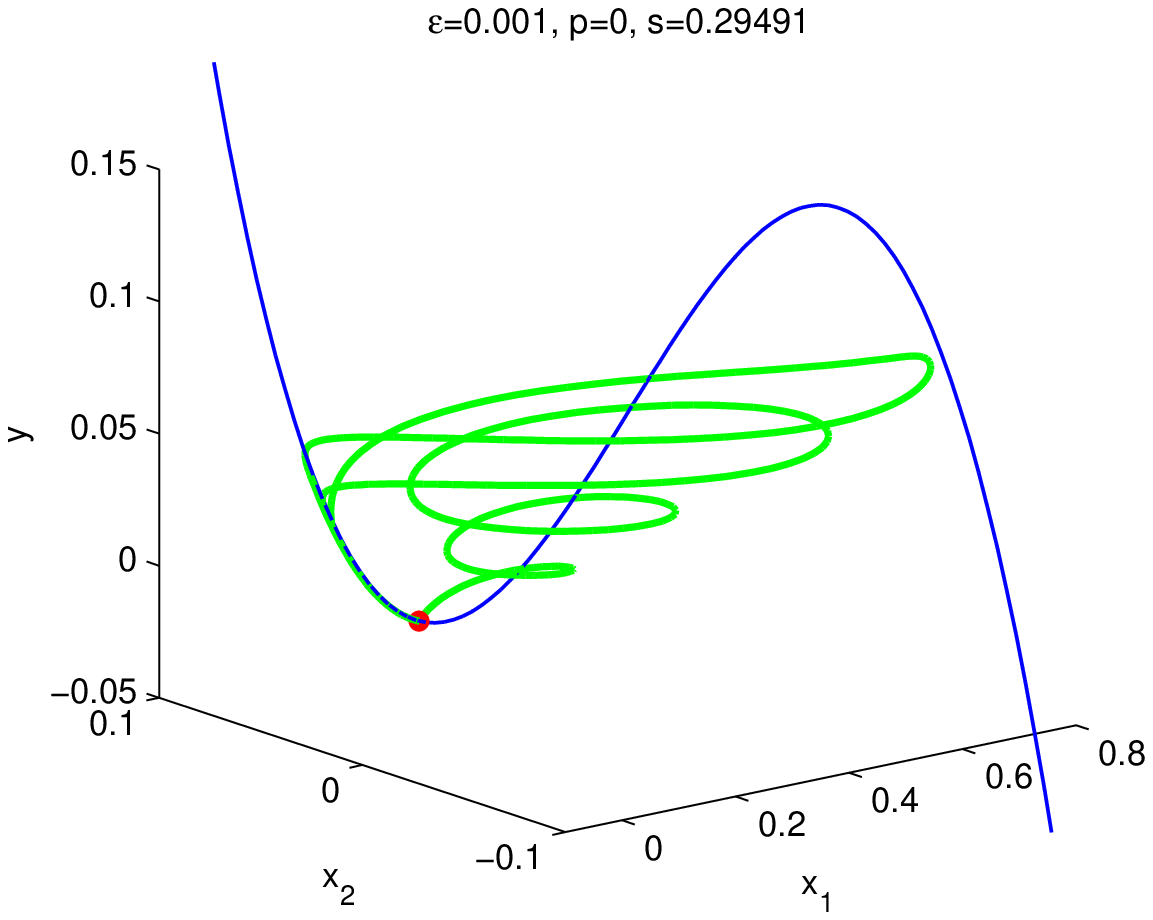}\label{fig:slow_wave2}}
\caption{Homoclinic orbits (green) representing slow waves in the FitzHugh-Nagumo equation. The slow manifold $S$ is shown in blue and the equilibrium $q$ in red.}
\label{fig:slow_waves}
\end{figure}

Now we consider the case of ``slow waves'' and work with smaller wave speeds $s$. Homoclinic orbits representing slow waves should be thought of as perturbations of singular limit orbits for the FitzHugh-Nagumo equation (\ref{eq:fhn}) with $s=0$. In this case the fast subsystem
\begin{eqnarray}
x_1'&=& x_2\nonumber\\
x_2'&=&\frac15 (-f(x_1)+y-p)
\end{eqnarray}
is Hamiltonian. Singular homoclinic orbits exist in a single fast subsystem with the y-coordinate of the equilibrium $y=x_1^*$. A direct application of Fenichel theory implies that a perturbed singular ``slow'' homoclinic orbit persists for $\epsilon>0$ \cite{Szmolyan1}. Again it is possible to compute parameter values $(p_1,s_1)$ at which homoclinic orbits for $\epsilon>0$ exist \cite{GuckCK1}. To compute the orbits themselves a similar approach as described above can be used. We have to track when $W^u(q)$ enters a small neighbourhood of $W^s(S_{l,\epsilon})$ respectively of $S_{l,\epsilon}$. Figure \ref{fig:slow_waves} shows two computed homoclinic orbits for $p_1=0$ and $s_1\approx 0.29491$.\\

The orbits spiral around the middle branch and do not enter the vicinity of $S_{r,\epsilon}$. This is expected as the middle branch $S_{m}$ of the critical manifold consists of unstable spiral equilibria for the fast subsystems. The Hamiltonian analysis for the case $s=0$ shows that the singular slow homoclinic orbits are not close to an equilibrium on $S_r$ for values of $p$ approximately between $-0.24$ and $0.05$ (see \cite{GuckCK1}). In Figure \ref{fig:slow_wave1} a homoclinic orbit enters the vicinity of the slow manifold $S_{l,\epsilon}$ and returns directly to $q$. Figure \ref{fig:slow_wave2} shows a homoclinic orbit that makes one additional large excursion around $S_{m,\epsilon}$ after it was close to $S_{r,\epsilon}$ and then returns to $q$; hence we refer to the orbit in \ref{fig:slow_wave2} as a double pulse homoclinic orbit. The same double pulse phenomenon exists for fast waves as well. In this case the double pulse orbit has no additional interaction with the middle branch $S_m$ and therefore it is difficult to distinguish between different pulse types for fast waves numerically and graphically as the second loop follows the first one very closely.  

\subsection{A Model of Reciprocal Inhibition}

This example demonstrates the use of our algorithm to compute trajectories in 
saddle-type slow manifolds of systems with two slow variables. The model is a 
caricature of a pair of neurons that are coupled with \emph{reciprocal inhibition}~\cite{rowat}.
The vector field is
\begin{equation}
\begin{split}
v_1' & = - \left(v_1 - a \tanh\left( \frac{\sigma_1
v_1}{a}  \right) + q_1 + \omega f(v_2)  (v_1-r) \right) \\
v_2' & = - \left(v_2  - a \tanh\left( \frac{\sigma_2
v_2}{a}   \right)+ q_2 + \omega f(v_1)  (v_2-r) \right) \\
q_1' & = \eps (-q_1 + s v_1) \\
q_2' & = \eps (-q_2 + s v_2)\\
f(x) & =\frac{1.0}{1.0+ \exp (-4 \gamma (x-\theta )) }
\end{split}
\label{rieqn}
\end{equation}
In this model, $v_1$ and $v_2$ are interpreted as the membrane potential of two neurons that are coupled synaptically through the terms involving $f$. The variables $q_1$ and $q_2$ represent the gating of membrane channels in the neurons. The model is a caricature in that it does not incorporate the fast membrane currents which give rise to action potentials. Still more reduced models~\cite{wang,skinner} have been used to study reciprocal inhibition of a pair of neurons. Reciprocal inhibition between a pair of identical neurons has long been viewed as a mechanism for generating repetitive alternating activity in motor systems~\cite{Brown}. Guckenheimer, Hoffman and Weckesser~\cite{GHW} investigated the properties of this model when the two neurons have different parameters and therefore are not identical.
They observed that canards of several kinds were encountered while continuing periodic orbits with AUTO calculations. The bifurcation mechanisms encountered
in these continuation studies still have not been identified despite intensive
efforts by Lust~\cite{Lust} to compute the multipliers of periodic orbits accurately.
Our algorithm for computing invariant slow manifolds of saddle-type provides a promising new tool for investigating the bifurcations that take place in this system with two slow and two fast variables. Here we illustrate that the algorithm is indeed 
capable of computing trajectories that lie on these manifolds without pursuing bifurcation analysis of the system.

\begin{figure}[!htb]
  \begin{center}
\includegraphics[height=3in]{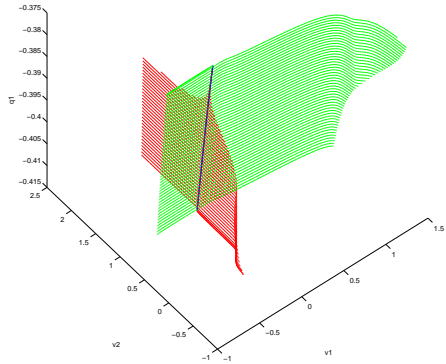}
 \end{center}
  \caption{A trajectory $\gamma$ on a saddle-type slow manifold of system \eqref{rieqn}. The trajectory is drawn black. Trajectories with initial conditions that are displaced by $10^{-8}$ from $\gamma$ along the direction of its strong stable manifolds are drawn red, while trajectories with initial conditions that are displaced by $10^{-8}$ from $\gamma$ along the direction of its strong unstable manifolds are drawn green. }
\label{rifig}
\end{figure}

The periodic trajectory discussed in section 5.2 of \cite{GHW} has three different canard segments. Here we focus on the segment labelled $B$ in Figure 6(c) of \cite{GHW}. The segment is a \emph{fold-initiated} canard that begins as a fast
trajectory, flows near a fold of the critical manifold and then moves along a saddle-type sheet of the slow manifold. Both neurons in the model have parameter values $\omega= 0.03, \, \gamma = 10, \, r= -4, \, \theta =0.01333, \, a= 1, \, s= 1 $
while $\sigma_1 = 3$ and $\sigma_2 = 1.2652372051$. One of the points $p$ on the segment $B$ has coordinates
$(-0.16851015831,0.85854544475,-0.41290838536,-0.062963871)$. We projected $p$ onto the critical manifold along the $q$ directions retaining the $v$ coordinates of $p$
and computed a trajectory $\gamma_{slow}$ of the slow flow on the critical manifold
with this initial condition. While the slow flow is an algebraic-differential equation, the critical manifold of \eqref{rieqn} is easily written as a graph of a function $q = h(v)$ and the slow flow equations can be written as a vector field in $v$. The trajectory $\gamma_{slow}$ was taken as input for our algorithm. Boundary conditions were selected so that the initial point of the trajectory $\gamma$ retains the same $v$ coordinates as $p$. Figure~\ref{rifig} displays the trajectory $\gamma$ obtained from our algorithm in black together with trajectories of its strong stable and unstable manifolds. The distance of the initial conditions for the trajectories on the strong stable and unstable manifolds from $\gamma$ is $10^{-8}$. Note that the 
first trajectories of the strong stable manifold at the bottom of the figure both flow down and to the right, reflecting that the initial points of these trajectories  do not straddle the slow manifold in the strong stable direction. Similar behavior occurs at the final point of $\gamma$ in the strong unstable direction. This behavior is to be expected because the boundary conditions constrain the strong stable coordinate of the first point of $\gamma$ to have a value close to that on the critical manifold rather than the invariant slow manifold. At the final point of $\gamma$, the strong unstable coordinate is determined by the critical manifold.
The behavior of $\gamma$ is what we expect from our algorithm: the computed trajectory approaches the slow manifold of saddle-type along a strong stable direction at its beginning, flows along the slow manifold to a high degree of accuracy to near its end and then leaves the slow manifold along a strong unstable direction. The length of $\gamma$ is much longer than the segment $B$ shown in Figure 6(c) of \cite{GHW}. 
 
{\bf Acknowledgment:} This research was partially supported by grants from the Department of Energy and the National Science Foundation.


\begin{thebibliography}{99}

\bibitem{Bowen} 

\bibitem{Brown} T.\ Brown (1911), The intrinsic factors in the act of progression
in the mammal, Proc.\ Roy.\ Soc.\ Lond.\ B 84:308-319.

\bibitem{Champ}
A.R. Champneys, V.~Kirk, E.~Knobloch, B.E. Oldeman, and J.~Sneyd (2007),
\newblock When shil'nikov meets hopf in excitable systems.
\newblock {\em SIAM Journal of Applied Dynamical Systems}, 6(4).

\bibitem{Fenichel} N. Fenichel (1971), Persistence and smoothness of invariant
manifolds for flows, Indiana Univ. Math. J. 21, 193-226.

\bibitem{FitzHugh}
R.~FitzHugh.
\newblock Mathematical models of threshold phenomena in the nerve membrane.
\newblock {\em Bull. Math. Biophysics}, 17:257--269, 1955.

\bibitem{Gorman}
A L Gorman and M V Thomas (1978) Changes in the intracellular concentration of free calcium ions in a pace-maker neurone, measured with the metallochromic indicator dye arsenazo III, J Physiol. 275:357-76.

\bibitem{GH} J. Guckenheimer and P. Holmes (1983), Nonlinear Oscillations, 
Dynamical Systems, and Bifurcation of Vector Fields, Springer Verlag.

\bibitem{GHW} Guckenheimer J, Hoffman K and  Weckesser W (2000), Numerical 
computation of canards, Int. J. Bif. Chaos 10, 2669--87 

\bibitem{GL}
J. Guckenheimer and D. LaMar,
Periodic orbit continuation in multiple time scale systems. Numerical continuation methods for dynamical systems, 253--267,
Underst. Complex Syst., Springer, Dordrecht, 2007.

\bibitem{GM}
J. Guckenheimer and B. Meloon (2000), Computing Periodic Orbits and their
Bifurcations with Automatic Differentiation, SIAM J. Sci. Comp.,
22, 951-985.

\bibitem{Montreal02}
J. Guckenheimer, Bifurcations of Relaxation Oscillations, in Normal Forms,
Bifurcations and Finiteness Problems in Differential Equations, Y. Ilyashenko
and C. Rousseau, eds. Kluwer, 295-316, 2004.

\bibitem{HW}
E. Hairer and G. Wanner, Solving ordinary differential equations. II. Stiff and differential-algebraic problems. Second edition. Springer Series in Computational Mathematics, 14. Springer-Verlag, Berlin, 1996.

\bibitem{Henon}
M. Hénon(1976), A two-dimensional mapping with a strange attractor.  Comm. Math. Phys.  50:69--77. 

\bibitem{Iz}
E. Izhikevich (2000), Neural Excitability, Spiking, and Bursting, Int. J. Bif. Chaos, 10:1171--1266.

\bibitem{Jones}
C. Jones,
Geometric singular perturbation theory. Dynamical systems (Montecatini Terme, 1994), 44--118,
Lecture Notes in Math., 1609, Springer, Berlin, 1995. 

\bibitem{JK}
C. Jones and N. Kopell (1994), Tracking invariant manifolds with differential forms in singularly perturbed systems.  J. Differential Equations  108:64--88. 

\bibitem{JonesKopellLanger}
C.~Jones, N.~Kopell, and R.~Langer.
\newblock Construction of the fitzhugh-nagumo pulse using differential forms.
\newblock {\em in: Multiple-Time-Scale Dynamical Systems}, pages 101--113,
  2001.

\bibitem{Ter1}
E. Lee and D. Terman (1999) Uniqueness and stability of periodic bursting solutions,
J. Diff. Eq., 158:48--78.

\bibitem{Lust}
K. Lust (2001), Improved numerical Floquet multipliers.  Internat. J. Bifur. Chaos Appl. Sci. Engrg.  11 , 2389--2410.

\bibitem{KS}
E. Kandel, J. Schwartz and T. Jessell (2000) Principles of Neuroscience, McGraw-Hill.

\bibitem{GuckCK1}
Christian Kuehn and John Guckenheimer.
\newblock Homoclinic orbits of the fitzhugh-nagumo equation: The singular
  limit.
\newblock {\em submitted}, 2008.

\bibitem{ML}
Morris C. and Lecar H. (1981), Voltage oscillations in the barnacle giant
muscle fiber. Biophysical Journal, 35: 193-213.

\bibitem{Nagumo}
J.~Nagumo, S.~Arimoto, and S.~Yoshizawa.
\newblock An active pulse transmission line simulating nerve axon.
\newblock {\em Proc. IRE}, 50:2061--2070, 1962.

\bibitem{Rinzel} J. Rinzel (1987), A formal classification of bursting 
mechanisms in excitable systems, Proc. Intern. Congr. of Mathematicians 
(A.M. Gleason, ed.), Amer. Math. Soc., 1578-1594.

\bibitem{RE} J. Rinzel and B. Ermentrout (1989), Analysis of neural 
excitability and oscillations, in Methods of Neural Modeling : From 
Synapses to Networks, C. Koch and I. Segev, eds., MIT Press, 135-169.

\bibitem{rowat}
Peter~F. Rowat and Allen~I. Selverston (1993),
\newblock Modeling the gastric mill central pattern generator of the lobster
  with a relaxation-oscillator network.
\newblock {\em Journal of Neurophysiology}, 70(3):1030--1053.

\bibitem{skinner}
Frances~K. Skinner, Nancy Kopell, and Eve Marder (1994),
\newblock Mechanisms for oscillation and frequency control in reciprocally
  inhibitory model neural networks.
\newblock {\em Journal of Computational Neuroscience}, 1:69--87.

\bibitem{Szmolyan1}
Peter Szmolyan.
\newblock Transversal heteroclinic and homoclinic orbits in singular
  perturbation problems.
\newblock {\em Journal of Differential Equations}, 92:252--281, 1991.

\bibitem{Ter}
D. Terman (1991), Chaotic spikes arising from a model of bursting in excitable membranes, SIAM J. Appl. Math. 51:1418--1450.

\bibitem{wang}
Xiao jing Wang and John Rinzel (1992),
\newblock Alternating and synchronous rhythms in reciprocally inhibitory model
  neurons.
\newblock {\em Neural Computation}, 4:84--97.


\end{thebibliography}
\end{document}